\documentclass[leqno]{article}
\usepackage{verbatim, amsmath, amscd, amssymb}
\usepackage[all]{xy}
\newcommand{\cI}{{\mathcal I}}

\newcommand{\cC}{{\mathcal C}}
\newcommand{\cD}{{\mathcal D}}
\newcommand{\cE}{{\mathcal E}}
\newcommand{\cF}{{\mathcal F}}

\newcommand{\cK}{{\mathcal K}}

\newcommand{\cM}{{\mathcal M}}

\newcommand{\cO}{{\mathcal O}}
\newcommand{\cP}{{\mathcal P}}
\newcommand{\cQ}{{\mathcal Q}}

\newcommand{\cV}{{\mathcal V}}
\newcommand{\cZ}{{\mathcal Z}}

\newcommand{\fm}{\mathfrak m}
\newcommand{\frn}{\mathfrak n}

\newcommand{\rd}{\mathrm{d}}

\newcommand{\rT}{\mathrm{T}}

\newcommand{\bbC}{\mathbb C}

\newcommand{\bbP}{\mathbb P}
\newcommand{\bbN}{\mathbb N}
\newcommand{\bbQ}{\mathbb Q}
\newcommand{\bbR}{\mathbb R}

\newcommand{\bbZ}{\mathbb Z}

\newcommand{\Ad}{\mathrm{Ad}}
\newcommand{\Ann}{\mathrm{Ann}}

\newcommand{\ch}{\mathrm{ch}\,}

\newcommand{\Dist}{\mathrm{Dist}}

\newcommand{\Ext}{\mathrm{Ext}}

\newcommand{\Frac}{\mathrm{Frac}}

\newcommand{\hd}{\mathrm{hd}}

\newcommand{\id}{\mathrm{id}}
\newcommand{\im}{\mathrm{im}}

\newcommand{\ind}{\mathrm{ind}}
\newcommand{\coind}{\mathrm{coind}}

\newcommand{\Lie}{\mathrm{Lie}}

\newcommand{\bfmodgr}{\mathbf{modgr}}

\newcommand{\op}{\mathrm{op}}

\newcommand{\pr}{\mathrm{pr}}

\newcommand{\rad}{\mathrm{rad}}

\newcommand{\Ru}{\mathrm{Ru}}
\newcommand{\soc}{\mathrm{soc}}

\newcommand{\Mod}{\mathbf{Mod}}
\newcommand{\Sch}{\mathbf{Sch}}

\newcommand{\lbr}{\begin{bmatrix}}
\newcommand{\rbr}{\end{bmatrix}}
\newcommand{\for}{\bigcirc\kern-2.6ex \because}
\newcommand{\forb}{\bigcirc\kern-2.8ex \because}
\newcommand{\forbb}{\bigcirc\kern-3.0ex \because}
\newcommand{\forbbb}{\bigcirc\kern-3.1ex \because}

\newcommand\pf{\noindent {\bf Proof:  }}
\newtheorem{thm}{Theorem:}

\newtheorem{prop}{Proposition:}

\newtheorem{lem}{Lemma:}

\newtheorem{cor}{Corollary:}

\begin{document}
\large
\title{
{\bf
Loewy series of
\\
parabolically induced $G_1T$-Verma modules
}
\thanks{supported in part by JSPS Grants in Aid for Scientific Research 23740004 AN and 23540023 KM}
\author{
A\textsc{be} Noriyuki
\\
Hokkaido University
\\
Creative Research Institution (CRIS)
\\
abenori@math.sci.hokudai.ac.jp
\\
\and
K\textsc{aneda} Masaharu
\\
Osaka City University
\\
Department of Mathematics
\\
kaneda@sci.osaka-cu.ac.jp
}
}
\maketitle

\begin{abstract}
Assuming the Lusztig conjecture on the irreducible characters for reductive algebraic groups in positive characteristic
$p$, 
which is now a theorem for large $p$,
we show that the modules for their Frobenius kernels induced from the simple modules of $p$-regular highest weights for their parabolic subgroups are rigid and determine their Loewy series. 
\end{abstract}

Let $G$ be a reductive algebraic group over an algebraically closed field 
$\Bbbk$ of positive characteristic
$p$, $P$ a parabolic subgroup of $G$, $T$ a maximal torus of $P$,
and
$G_1$
(resp. $P_1$) the Frobenius kernel of
$G$
(resp. $P$).
In this paper we study the structure of $G_1T$-modules
induced from the simple $P_1T$-modules of $p$-regular highest weights.
Thus
our study goes parallel to
parabolically induced
Verma modules in characteristic $0$.
In case $P$ is a Borel subgroup of $G$, assuming Lusztig's conjecture for the irreducible characters for $G_1T$,
which is now a theorem for large $p$
thanks to
\cite{AJS},
\cite{KL},
\cite{L94}, \cite{KT},
or more recently to
\cite{F},
H. H. Andersen and the second author of the present paper
showed 
that the induced modules are rigid and determined
their Loewy series
\cite{AK89}.
We now show that the parabolically induced modules
are also rigid and
describe
their Loewy series.

To go into more details, let $B$ be a Borel subgroup of $P$
containing
$T$,
$\Lambda$
the character group of $B$, $R\subset\Lambda$  
the root system of $G$ relative to
$T$, and 
$R^+$ the positive system of
$R$ such that the roots of $B$ are
$-R^+$.
We let $R^s$ denote the set of simple roots,
and
$I$ a subset of
$R^s$  
such that
the root subgroups
$U_\alpha$ 
of $G$ associated to 
$\alpha\in I$
generate $P$ together with
$B$.
Denote by
$\hat\nabla_P
$ the induction functor from the category of $P_1T$-modules to the category of
$G_1T$-modules, and
let
$\hat L^P(\lambda)$
denote the simple $P_1T$-module of highest weight $\lambda\in \Lambda$.
Our object of study is
$\hat\nabla_P(\hat L^P(\lambda))$.
After stating some generalities in \S\S 1 and 2,
we specialize into the case where
$\lambda$ is $p$-regular, i.e., if $\alpha^\vee$ is the coroot of
each root
$\alpha$
and if
$\rho=\frac{1}{2}\sum_{\alpha\in R^+}\alpha$, the case when
$p\not|\langle\lambda+\rho,
\alpha^\vee\rangle$ for all 
roots $\alpha$.
If $M$ is a finite dimensional
$G_1T$-module,
we call the sum of its simple submodules the socle of
$M$ and denote it by
$\soc M=\soc^1M$.
If $\pi:M\to
M/\soc M$ is the quotient,
we let
$\soc^2M=\pi^{-1}\soc(M/\soc M)$, and repeat to construct a filtration
$0<\soc M<\soc^2M<\dots <M$,
called the socle series of $M$.
Dually,
we call
the intersection of all its maximal submodules the radical of $M$ and denote it
by
$\rad M=\rad^1M$. Letting
$\rad^iM=\rad(\rad^{i-1}M)$ for
$i>1$, one obtains a filtration
$M>\rad M>\rad^2M>\dots>0$,
called the radical series of
$M$.
It is known that the minimal $i$
such that $\soc^iM=M$
and the minimal
$j$ such that
$\rad^jM=0$ coincide, called the Loewy length of $M$ and denoted
$\ell\ell(M)$.
By definition each $\soc_iM=\soc^iM/\soc^{i-1}M$,
called the $i$-th socle layer
of $M$,
and $\rad_iM=\rad^iM/\rad^{i+1}M$,
called the $i$-th radical layer of $M$, are semisimple.
Any filtration $0<M^1<M^2<\dots<M$
with each subquotient semisimple has the length
at least
$\ell\ell(M)$.
If the length of the filtration such a filtration $M^\bullet$ is $\ell\ell(M)$, then
$\soc^iM\geq M^i\geq\rad^{\ell\ell(M)-i+1}M$ for each $i$.
We say $M$ is rigid iff
the socle series and the radical series of $M$ coincide.
In \S3 we employ
graded representation theory from
\cite{AJS}
to show that
the induction functor
$\hat\nabla_P$ 
is
$\bbZ$-graded.
Each block of finite dimensional
$G_1T$-modules 
is equivalent to the
category of $p\bbZ R$-graded 
finite dimensional modules
over the endomorphism $\Bbbk$-algebra
$E$ of
a projective
$p\bbZ R$-generator of the block.
Assuming Lusztig's conjecture for the irreducible characters for $G_1T$,
Andersen, Jantzen and Soergel
\cite{AJS} showed that
the algebra $E$ for
a $p$-regular block
is $(p\bbZ R\times\bbZ)$-graded
and is Koszul with respect to its
$\bbZ$-gradation.
We 
show in \S4 that
the rigidity of $\hat\nabla_P(\hat L^P(\lambda))$ 
for $p$-regular
$\lambda$
follows from 
a result in
\cite{BGS}. 
Unlike the case $P=B$ the number of $G_1T$-composition factors of
$\hat\nabla_P(\hat L^P(\lambda))$ varies depending on
the highest weight
$\lambda$.
Nonetheless, 
we show also in \S4 that the Loewy
length of
$\hat\nabla_P(\hat L^P(\lambda))$ is uniformly
$\ell(w_0w_I)+1$
with
$w_0$
(resp. $w_I$)
the longest element of the Weyl group
$W$
(resp. $W_I$) of $G$
(resp. $P$).
In \S5 we determine the Loewy series of 
$\hat\nabla_P(\hat L^P(\lambda))$.

Given a category $C$ and its objects $X$ and $Y$, $C(X, Y)$ will denote the set of morphisms in
$C$ from $X$ to $Y$.

The second author of the paper is grateful to Arun Ram and the University of Melbourne for the hospitality during his visit in August of 2011, where some of initial ideas of the present work were conceived.
We also thank Peter Fiebig for a helpful discussion.

\begin{center}
$1^\circ$
{\bf Some generalities}
\end{center} 

Let 
$G$ be a reductive algebraic group
over an algebraically field
$\Bbbk$ of positive characteristic
$p$, $B$ a Borel subgroup of $G$,
$T$ a maximal torus of $B$,
$\Lambda$ 
the character group of $B$, $R\subset\Lambda$ 
the root system of $G$ relative to
$T$, and 
$R^+$ 
the positive system of
$R$ such that the roots of $B$ are
$-R^+$.
We let $R^s\subset R^+$ denote the set of simple roots,
and
$\Lambda^+$ 
the set of dominant weights of
$\Lambda$.
For each $\alpha\in R$ we let 
$\alpha^\vee$ denote the coroot of
$\alpha$.
Let
$W$ be the Weyl group of $G$
generated by the reflections
$s_\alpha$, $\alpha\in R$,
and $\ell$ 
the length function on $W$ with respect to the simple reflections.
Let $w_0$ be the longest element of $W$.

For each $\alpha\in R$ let $U_\alpha$ denote the root subgroup of $G$ associated to
$\alpha$.
Let
$I\subseteq R^s$ and
$P=P_I=\langle
B, U_\alpha\mid\alpha\in I\rangle$ 
the standard parabolic subgroup of $G$
associated to
$I$, and let $L_I$ denote its standard Levi subgroup.
Let
$R_I\subseteq R$ denote the root system of $L_I$
with its induced positive system
$R^+_I$.
Put $\Lambda_P=
\{\lambda\in
\Lambda\mid
\langle\lambda,\alpha^\vee\rangle=0
\ \forall\alpha\in I\}$
and
$\Lambda_I^+=\{\lambda\in
\Lambda\mid
\langle\lambda,\alpha^\vee\rangle\geq0
\ \forall\alpha\in I\}$.
Let $W_I$ be the Weyl group of $P$ and $w_I$ 
its longest element.
Put $w^I=w_0w_I$.
Let
$\rho=\frac{1}{2}\sum_{\alpha\in R^+}\alpha$
and
$\rho_P=\frac{1}{2}\sum_{\alpha\in
R^+\setminus R_I^+}\alpha\in\Lambda\otimes_\bbZ\bbQ$.
For simplicity we will assume
$G$ is semisimple and simply connected.
Let
$W_{p}=W\ltimes
p\bbZ
R$,
$W_{I,p}=W_I\ltimes
p\bbZ
R_I$,
and $\rho_I=\frac{1}{2}\sum_{\alpha\in R^+_I}\alpha=\rho-\rho_P$. 
For $x\in
W_p$  we will write
$x\bullet\lambda$
for
$x(\lambda+\rho)-\rho$.
In case $x\in W_{I,p}$,
$x\bullet\lambda=x(\lambda+\rho_I)-\rho_I$.
We will also let
$(-x)\bullet\lambda$ stand for
$-(x\bullet\lambda)-2\rho=-x(\lambda+\rho)-\rho$.

\noindent
(1.1)
Let $\alpha_0$ be the highest short root of
$R$ and let
 $h=\langle\rho,\alpha_0^\vee\rangle+1$ the Coxeter number of
$G$.

\begin{lem}
$2\rho_P=w_I\rho+\rho=w_0(w^I\bullet0)\in\Lambda_P\cap\Lambda^+$ with
$\langle2\rho_P,\alpha^\vee\rangle\in[2,h]$
$\forall\alpha\in R^s\setminus
I$.

\end{lem}

\pf
One has
\begin{align*}
w_0(w^I\bullet0)
&=
w_0(w_0w_I\rho-\rho)
=
w_I\rho+\rho
=
\rho+w_I\frac{1}{2}(\sum_{\beta\in R^+\setminus R_I^+}\beta
+
\sum_{\beta\in R_I^+}\beta)
\\
&=
\frac{1}{2}(\sum_{\beta\in R^+\setminus R_I^+}\beta
+
\sum_{\beta\in R_I^+}\beta)+
\frac{1}{2}(\sum_{\beta\in R^+\setminus R_I^+}\beta
-
\sum_{\beta\in R_I^+}\beta)
=
\sum_{\beta\in R^+\setminus R_I^+}\beta
=2\rho_P.
\end{align*}

If $\alpha\in I$,
$\langle2\rho_P,\alpha^\vee\rangle=
\langle
w_I\rho+\rho,\alpha^\vee\rangle
=
\langle
\rho,w_I\alpha^\vee\rangle+
1=0$,
and hence
$2\rho_P\in\Lambda_P$.
If
$\alpha\in R^s\setminus
I$,
$\langle2\rho_P,\alpha^\vee\rangle
=
\langle
w_I\rho+\rho,\alpha^\vee\rangle
=
\langle
\rho,w_I\alpha^\vee\rangle+1
\leq
\langle\rho,\alpha_0^\vee\rangle+1
=
h$.

\noindent
(1.2)
If $H\leq K$ are closed
subgroups of $G$,
we let $\ind_H^K$ denote the induction functor from the category 
$H\Mod$
of
rational $H$-modules
to the category $K\Mod$
of rational $K$-modules:
if $M\in H\Mod$,
$\ind_H^KM=
\{f\in\Sch_\Bbbk(K,M)\mid
f(kh)=h^{-1}f(k)\ \forall k\in K
\forall h\in H\}$.
We let $\Dist(H)$
(resp. $\Dist(K)$)
denote the algebra
of distributions on $H$
(resp. $K$)
and let
$\coind_H^K=\Dist(K)\otimes_{\Dist(H)}?$
denote the coinduction functor from
$\Dist(H)\Mod$ to $\Dist(K)\Mod$. 
For a finite dimensional $H$-module $M$
 we will mean by $M^*$
 the
$\Bbbk$-linear dual of
$M$.
By $\otimes$ we will always mean
$\otimes_\Bbbk$ 
unless otherwise specified.
Let $H_1$ denote the Frobenius kernel of
$H$.

If $M$ is a $P$-module,
$\coind_{P_1}^{G_1}M$ extends to
a $G_1P$-module with
$P$ acting on
$\Dist(G_1)$ and $\Dist(P_1)$
by the adjoint action
and as given on
$M$,
in which case we will write
$\coind_P^{G_1P}M$ for
$\coind_{P_1}^{G_1}M$
\cite[I.8.20]{J}.
Let
$\Ru(P)$ 
denote
the unipotent radical of
$P$.

\noindent
{\bf Proposition
(cf. \cite[II.3.5]{J}):}
{\it 
Let
$M\in P\Mod$.

(i)
There is an isomorphism of
$G_1P$-modules
$\ind_P^{G_1P}M\simeq
\coind_{P}^{G_1P}(M\otimes
2(1-p)\rho_P)$.

(ii)
If $M$ is finite dimensional,
there is an isomorphism of
$G_1P$-modules
\[
(\ind_P^{G_1P}M)^*
\simeq
\ind_{P}^{G_1P}(M^*\otimes
2(p-1)\rho_P).
\]

}

\pf
Recall from \cite[I.8.20]{J} an isomorphism of
$G_1P$-modules
\begin{align}
\coind_{P}^{G_1P}M
&\simeq
\ind_P^{G_1P}(M\otimes\chi\vert_P(\chi')^{-1}),
\\
(\ind_{P}^{G_1P}M)^*
&\simeq
\ind_P^{G_1P}(M^*\otimes\chi\vert_P(\chi')^{-1})
\quad\text{if $\dim M<\infty$},\end{align}
where
$\chi$
(resp. $\chi'$)
is a 1-dimensional representation of
$G$
(resp. $P$) through which 
$G$
(resp. $P$)
acts on 
$\Dist(G_1)^{G_1}_\ell
=\{\mu\in\Dist(G_1)|\rho_\ell(x)\mu=\mu\ \forall x\in G_1\}
$
(resp.
$\Dist(P_1)^{P_1}_\ell=\{\mu\in\Dist(P_1)|\rho_\ell(x)\mu=\mu\ \forall x\in P_1\}
$),
where
$\rho_\ell$ denotes the left regular action.
As $\chi$ is trivial by \cite[II.3.4/I.9.7]{J},
(1) and (2) read, resp.,
\begin{align}
\coind_{P}^{G_1P}M
&\simeq
\ind_P^{G_1P}(M\otimes(\chi')^{-1}),
\\
(\ind_{P}^{G_1P}M)^*
&\simeq
\ind_P^{G_1P}(M^*\otimes(\chi')^{-1}).
\end{align}

Recall from
\cite[I.9.7]{J} that
$\chi'$
is given by
$g\mapsto
\det(\Ad(g))^{p-1}$,
$g\in P$.
In particular, $\chi'$ factors through
$P/\Ru(P)$, and is trivial on the 
derived subgroup of
$L_I$.
To compute
$\chi'$, therefore, we have only to consider
the adjoint representation of $T$ on  $\Lie(P)=
\Lie(T)\oplus
\bigoplus_{\beta\in R^+}\Lie(U_{-\beta})
\oplus
\bigoplus_{\alpha\in R^+_I}\Lie(U_\alpha)$.
Thus for each $t\in T$
\begin{align*}
\det(\Ad(t))
&=
(\sum_{\beta\in R^+}-\beta+\sum_{\alpha\in R^+_I}\alpha)(t)
=
\{-(\sum_{\beta\in R^+}\beta-\sum_{\alpha\in R^+_I}\alpha)\}(t)
=
(-\sum_{\beta\in R^+\setminus
R^+_I}\beta)(t)
\\
&=
(-2\rho_P)(t).
\end{align*}
It follows that
$\chi'=(p-1)(-2\rho_P)$, and hence the assertions.

\setcounter{equation}{0}
\noindent
(1.3)
Likewise,
write 
$\coind_{P_1T}^{G_1T}M$
for the $G_1T$-module
$\coind_{P_1}^{G_1}M$
in case $M$ is a $P_1T$-module.

\begin{prop}
Let
$M\in P_1T\Mod$.

(i)
There is an isomorphism of
$G_1T$-modules
$\ind_{P_1T}^{G_1T}M\simeq
\coind_{P_1T}^{G_1T}(M\otimes
2(1-p)\rho_P)$.

(ii)
If $M$ is finite dimensional,
there is an isomorphism of
$G_1T$-modules
\[
(\ind_{P_1T}^{G_1T}M)^*
\simeq
\ind_{P_1T}^{G_1T}(M^*\otimes
2(p-1)\rho_P).
\]

\end{prop}

\noindent
(1.4)
If $L$ is a simple $P$-module,
the $P$-action on $L$ factors through 
$P/\Ru(P)$,
affording a simple $L_I$-module of highest weight belonging to
$\Lambda_I^+$.
For each $\lambda\in\Lambda^+$
(resp.
$\lambda\in\Lambda_I^+$),
we let
$L(\lambda)$ 
(resp. $L^P(\lambda)$)
denote the simple $G$-
(resp. $P$-)
module of
highest weight $\lambda$.
Likewise for
simple $P_1T$-modules.
For each $\lambda\in\Lambda$
we let
$\hat L(\lambda)$
(resp. $\hat L^P(\lambda)$)
denote the simple
$G_1T$-
(resp. $P_1T$-)
module of highest weight $\lambda$.
Let $\Lambda_p=\{\lambda\in\Lambda\mid
\langle\lambda,
\alpha^\vee\rangle\in[0,p[\ \forall\alpha\in R^s\}$.
Each
$\lambda\in\Lambda$
admits a decomposition
$\lambda=\lambda^0+p\lambda^1$
with
$\lambda^0\in\Lambda_p$ and
$\lambda^1\in\Lambda$.
Thus
$\hat L^P(\lambda)\simeq
L^P(\lambda^0)\otimes
p\lambda^1$; if $\lambda^0=\lambda^0_I+\lambda^c_I$ with $\lambda_I^c\in \Lambda_P$, then $\hat L^P(\lambda)\simeq L^P(\lambda^0_I)\otimes(\lambda_I^c+p\lambda^1)\simeq L^P(\lambda^0)\otimes p\lambda^1$.
In particular,
\begin{align}
\{
\ind_{P_1T}^{G_1T}
(\hat L^P(\lambda))
\}^*
\simeq
\ind_{P_1T}^{G_1T}(\hat
L^P((-w_I)\bullet\lambda)))\otimes
p(2\rho_P+w_I\lambda^1-\lambda^1)
\end{align}
with
$2\rho_P+w_I\lambda^1-\lambda^1\in\bbZ R$.

If $H$ is a closed subgroup of $G$ and if $M$ is an $H$-module,
we let
$\soc_HM$
(resp. $\rad_HM$)
denote the socle
(resp. the radical)
of
$M$, and put
$\hd_HM=M/(\rad_HM)$.

\begin{prop}
For each $\lambda\in \Lambda$
\begin{align*}
\soc_{G_1T}\ind_{P_1T}^{G_1T}(\hat L^P(\lambda))
&=
\hat L(\lambda),
\\
\hd_{G_1T}\ind_{P_1T}^{G_1T}(\hat L^P(\lambda))
&=
\hat L(-w_I\lambda^0-p\lambda^1+2(p-1)\rho_P)^*
\\
&\hspace{-2cm}
=
\hat L(w^I\bullet\lambda)\otimes
p(\lambda^1-2\rho_P-w_I\lambda^1+w_0((-w_I)\bullet\lambda)^1-((-w_I)\bullet\lambda)^1)).
\end{align*}

\end{prop}

\pf
For each
$\lambda\in\Lambda$
we have
$\soc_{P_1T}(\ind_{B_1T}^{P_1T}\lambda)\simeq
\soc_{P_1T}(\ind_{(B/\Ru(P))_1T}^{(P/\Ru(P))_1T}\lambda)=\hat L^P(\lambda)$.
Then
\[
\ind_{P_1T}^{G_1T}\hat L^P(\lambda)\leq
\ind_{P_1T}^{G_1T}\ind_{B_1T}^{P_1T}(\lambda)
\simeq
\ind_{B_1T}^{G_1T}\lambda.
\]
It follows that
$\soc_{G_1T}(\ind_{P_1T}^{G_1T}\hat L^P(\lambda))
=
\{\ind_{P_1T}^{G_1T}\hat L^P(\lambda)\}\cap
\soc_{G_1T}(
\ind_{B_1T}^{G_1T}\lambda)
=
\hat L(\lambda)$.
Then
\begin{align*}
\hd_{G_1T}(\ind_{P_1T}^{G_1T}\hat L^P(\lambda))
&\simeq
\{\soc_{G_1T}((\ind_{P_1T}^{G_1T}\hat L^P(\lambda))^*)\}^*
\\
&\simeq
\{\soc_{G_1T}(
\ind_{B_1T}^{G_1T}(\hat L^P(\lambda)^*\otimes
2(p-1)\rho_P))\}^*
\quad\text{by (1.2.ii)}.
\end{align*}
Now
$\hat L^P(\lambda)^*
=
(L^P(\lambda^0)\otimes
p\lambda^1)^*
=
L^P(\lambda^0)^*\otimes-p\lambda^1
=
L^P(-w_I\lambda^0)\otimes-p\lambda^1
=
\hat L^P(-w_I\lambda^0-p\lambda^1)$.
Also $\forall\nu\in\Lambda$
\begin{align*}
\hat L^P(\nu)\otimes2(p-1)\rho_P
&\leq
(\ind_{B_1T}^{P_1T}\nu)\otimes
2(p-1)\rho_P
\\
&\simeq
\ind_{B_1T}^{P_1T}(\nu\otimes
2(p-1)\rho_P)
\quad\text{by the tensor identity},
\end{align*}
and hence
\[
\hat L^P(\nu)\otimes2(p-1)\rho_P
\leq
\soc_{P_1T}
\ind_{B_1T}^{P_1T}(\nu\otimes
2(p-1)\rho_P)
=
\hat L^P(\nu\otimes
2(p-1)\rho_P).
\]
It follows that
\begin{align*}
\hd_{G_1T}(\ind_{P_1T}^{G_1T}\hat L^P(\lambda))
&\simeq
\{\soc_{G_1T}(
\ind_{B_1T}^{G_1T}(\hat L^P(-w_I\lambda^0-p\lambda^1+2(p-1)\rho_P)))\}^*
\\
&=
\hat L(-w_I\lambda^0-p\lambda^1+2(p-1)\rho_P)^*.
\end{align*}
Finally,
\begin{align*}
-w_I\lambda^0-
&
p\lambda^1+2(p-1)\rho_P
=
-w_I\lambda^0-p\lambda^1+(p-1)(w_I\rho+\rho)
\quad\text{by (1.1)}
\\
&=
-w_I(\lambda^0+\rho)-\rho+p(w_I\rho+\rho-\lambda^1)
=
(-w_I)\bullet\lambda+p(w_I\lambda^1+w_I\rho+\rho-\lambda^1)
\\
&=
(-w_I)\bullet\lambda+p(w_I\lambda^1+2\rho_P-\lambda^1)
\quad\text{by (1.1) again}.
\end{align*}
Thus
\begin{multline*}
\hat L(-w_I\lambda^0-
p\lambda^1+2(p-1)\rho_P)^*
=
\{
\hat L((-w_I)\bullet\lambda+
p(w_I\lambda^1+2\rho_P-\lambda^1))
\}^*
\\
=
\hat L(-w_0((-w_I)\bullet\lambda))\otimes
p\{-w_I\lambda^1-2\rho_P+\lambda^1+w_0((-w_I)\bullet\lambda)^1-
((-w_I)\bullet\lambda)^1\}
\end{multline*}
with
$
-w_0((-w_I)\bullet\lambda)=
-w_0(-w_I(\lambda+\rho)-\rho)
=
w_0w_I(\lambda+\rho)-\rho
=
w_0w_I\bullet\lambda=w^I\bullet\lambda$.

\noindent
(1.5)
{\bf Corollary:}
{\it   
Let $\lambda\in\Lambda$.

(i)
$\hd_{P_1T}\hat L(\lambda)=\hat L^P(\lambda)$
while
$\soc_{P_1T}\hat L(\lambda)=
\hat L^P(w_Iw_0\lambda^0+p\lambda^1)$.

(ii)
If $\lambda\in\Lambda^+$,
$\hd_{P}
L(\lambda)=L^P(\lambda)$
while
$\soc_{P}
 L(\lambda)=
L^P(w_Iw_0\lambda)$.

}

\pf
(i)
For each $\nu\in\Lambda$
\begin{align*}
P_1T\Mod(\hat L(\lambda),
\hat L^P(\nu))
&\simeq
G_1T\Mod(\hat L(\lambda),
\ind_{P_1T}^{G_1T}\hat L^P(\nu))
\\
&=
\delta_{\lambda\nu}\Bbbk
\quad\text{by (1.4)}.
\end{align*}
It follows that
$\hd_{P_1T}\hat L(\lambda)=\hat L^P(\lambda)$.
Then
\begin{align*}
\soc_{P_1T}\hat L(\lambda)
&\simeq
\{\hd_{P_1T}(\hat L(\lambda)^*)\}^*
=
\{\hd_{P_1T}\hat L(-w_0\lambda^0-p\lambda^1)\}^*
=
\hat L^P(-w_0\lambda^0-p\lambda^1)^*
\\
&=
\{
L^P(-w_0\lambda^0)\otimes-p\lambda^1\}^*
\simeq
L^P(w_Iw_0\lambda^0)\otimes
p\lambda^1
=
\hat L^P(w_Iw_0\lambda^0+p\lambda^1).
\end{align*}

(ii)
For each
$\mu\in\Lambda^+_I$
\begin{align*}
P\Mod(L(\lambda),
L^P(\mu))
&\simeq
G\Mod(L(\lambda),
\ind_{P}^{G}L^P(\mu))
\leq
G\Mod(L(\lambda),
\ind_{P}^{G}\ind_B^P(\mu))
\\
&\simeq
G\Mod(L(\lambda),
\ind_{B}^{G}(\mu))
=
\delta_{\lambda\mu}\Bbbk.
\end{align*}
It follows that
$\hd_PL(\lambda)=L^P(\lambda)$.
Then
\begin{align*}
\soc_PL(\lambda)
&\simeq
\{\hd_P(L(\lambda)^*)\}^*
=
\{\hd_PL(-w_0\lambda)
\}^*
=
L^P(-w_0\lambda)^*
=
L^P(w_Iw_0\lambda).
\end{align*}

\setcounter{equation}{0}
\noindent
(1.6)
Let $H$ be a closed subgroup of
$G$ and
$\phi$ an automorphism of
$H$.
If $M$ is an $H$-module, by ${^\phi \!M}$ we will mean an $H$-module of ambient
$\Bbbk$-linear space $M$ with the $H$-action twisted by
$\phi^{-1}$
\cite[I.2.15/3.5]{J}:
$\forall h\in H$
$\forall
m\in M$, the action of
$h$ on $m$ in ${^\phi \!M}$ is given by
$\phi^{-1}(h)m$.
In particular,
under the conjugate action of $W$ on $T$,
$\forall w\in W$ and $\forall\lambda\in \Lambda$,
\begin{equation}
^w\!\lambda=
w\lambda.
\end{equation}

If $K$ is a closed subgroup of $H$ and $V$ is a $K$-module,
there is an isomorphism of ${^w\!H}$-modules
\cite[I.3.5.4]{J}
\begin{equation}
^w
\ind_K^H(V)\simeq
\ind_{^wK}^{^wH}(^wV).
\end{equation}

Throughout the rest of the paper we will abbreviate
$\ind_{P_1T}^{G_1T}$
(resp.
$\ind_{B_1T}^{P_1T}$
)
as
$\hat\nabla_P$
(resp.
$\hat\nabla^P$).
More generally,
for $w\in W$ let
${^w\!P}=wPw^{-1}$
and put
$\hat\nabla_{^w\!P}=\ind_{(^w\!P)_1T}^{G_1T}$,
$\hat\nabla^{^w\!P}=\ind^{(^w\!P)_1T}_{(^w\!B)_1T}$.
Let also
$\hat\nabla_w=\ind_{(^wB)_1T}^{G_1T}$; we will abbreviate
$\hat\nabla_e$
as $\hat\nabla$.
For each
$\lambda\in\Lambda$
and
$w\in W$
we will write
$\lambda\langle w\rangle$
for
$\lambda+(p-1)(w\bullet0)$
after
\cite{AJS}.
Then
\begin{align}
^w
\hat\nabla_P(\hat L^P(\lambda))
&\simeq
\hat\nabla_{^w\!P}(^w\!\hat L^P(\lambda))
\quad\text{by (2)}
\\
\notag&\leq
\hat\nabla_{^w\!P}(^w
\hat\nabla^P(\lambda))
\\
\notag&\simeq
\hat\nabla_{^w\!P}(
\hat\nabla^{^w\!P}_{^w\!B}(^w
\lambda))
\quad\text{by (2) again}
\\
\notag&\simeq
\hat\nabla_{w}(
w
\lambda)
\quad\text{by (1)}
\\
\notag&=
\hat\nabla_{w}(
w\bullet
\lambda-w\bullet0)
\simeq
\hat\nabla_{w}(
w\bullet
\lambda+(p-1)(w\bullet0))\otimes
-p(w\bullet0)
\\
\notag
&=
\hat\nabla_{w}(
(w\bullet
\lambda)
\langle w\rangle)\otimes
-p(w\bullet0).
\end{align}

\setcounter{equation}{0}
\noindent
(1.7)
Put
$\hat\Delta=\coind_{B^+_1T}^{G_1T}$.
Let
$\tau$ be the 
Chevalley
antiinvolution of $G$ such that $\tau|_T=\id_T$
\cite[II.1.16]{J}, and hence
$\tau(U_\alpha)=U_{-\alpha}$
for each $\alpha\in R$.
If $H$ is a subgroup of $G$ and if $M$ is a finite dimensional $H$-module,
let
$M^\tau$ be the $\tau(H)$-module with the ambient space
$M^*$
and the
$\tau(H)$-action twisted by $\tau$:
$\forall x\in\tau(H)$,
$\forall f\in M^*$,
$\forall m\in M$,
$(xf)(m)=f(\tau(x)m)$.
Recall from
\cite[II.9.3.5]{J} that
there is a functorial isomorphism
$(?^\tau)\circ\hat\nabla\simeq
\hat \Delta\circ(?^\tau)$
on the category of
finite dimensional $B_1T$-modules.
More generally,
put
$B^+=\tau B$,
$P^+=\tau P=\langle B^+,U_{-\alpha}|\alpha\in I\rangle$
and let
$\hat\Delta_{P}=\coind_{P^+_1T}^{G_1T}$.
If
$M$ is a finite dimensional $P_1T$-module, there is an isomorphism of
$G_1T$-modules
\begin{equation}
(\hat\nabla_P(M))^\tau\simeq
\hat\Delta_{P}(M^\tau).
\end{equation}

Let
$U^+_1(w_I)=\prod_{\beta\in R^+\setminus R_I}U_{\beta,1}$
be the Frobenius kernel of
the unipotent radical of
$P^+$.
If $V$ is a $G_1T$-module,
let
$V^{U^+_1(w_I)}=
\{v\in V|xv=v\  \forall x\in
U^+_1(w_I)\}$.
If $M$ is a $B_1T$-module,
as
$G_1=U^+_1(w_I)P_1$,
$\hat\nabla(M)^{U^+_1(w_I)}
=
\{\Sch_\Bbbk(G_1T,M)^{B_1T}\}^{U^+_1(w_I)}
=
\Sch_\Bbbk(P_1T,M)^{B_1T}
$
maintains a structure of
$P_1T$-module such that
\begin{equation}
\hat\nabla(M)^{U^+_1(w_I)}=
\hat\nabla^P(M).
\end{equation}
Recall also that each 
$\hat\nabla(\lambda)$,
$\lambda\in\Lambda$,
is projective/injective as $B_1^+T$-module
\cite[II.9.5]{J}.
As
$U^+_1(w_I)$ is a normal subgroup of
$B^+_1$,
$U^+_1(w_I)$ is exact in $B^+_1$
\cite[I.6.5.2]{J}, and hence
$\hat\nabla(\lambda)$
remains injective/projective as
$U^+_1(w_I)$-module.

\begin{center}
$2^\circ$
{\bf Translation functors}
\end{center}

For
$\lambda,\mu\in \Lambda$
let
$\rT_\lambda^\mu$ denote the translation functor on 
the
$G_1T$-modules.
If $M$ is a $L_{I,1}T$-module,
we say
$M$ belongs to $\lambda$ iff
all the $L_{I,1}T$-composition factors of
$M$ are have highest weights belonging to
$W_{I,p}\bullet\lambda$.
We let $\rT_{I,\lambda}^\mu$
denote the translation functor on
the $L_{I,1}T$-modules.

For each $\alpha\in R$ and $n\in\bbZ$
let
$H_{\alpha,n}=\{v\in\Lambda \otimes_\bbZ\bbR\mid\langle
v+\rho,\alpha^\vee\rangle=pn\}$.
We call a connected component of
$\Lambda\otimes_\bbZ\bbR\setminus\cup_{\alpha\in R,n\in\bbZ}H_{\alpha,n}$
an alcove.
If $F\subseteq\Lambda\otimes_\bbZ\bbR$, $\overline F$ will
denote the closure of $F$ in
$\Lambda\otimes_\bbZ\bbR$.
We say
$\lambda\in \Lambda$ is $p$-regular iff
$\lambda$ lies in an alcove.
If $x\in W_p$ and $A$ is an alcove,
we will write
$xA$ to mean $x\bullet A$.

\noindent
(2.1)
{\bf Lemma:}
{\it
Let
$\eta\in\Lambda$ and
$E$ a simple $G$-module of extremal weight $\eta$.
If $w\eta\in\Lambda^+_I$,
$w\in W_I$, and if $\alpha\in I$,
then
$w\eta+\alpha$ is not a weight of
$E$.

}

\pf
Let
$x\in W$ with
$x\eta\in\Lambda^+$, and put
$\nu=x\eta$, $\nu'=w\eta$.
Let
$J=\{\beta\in I\mid
\langle\nu', \beta^\vee\rangle=0\}$,
$W_J=\langle s_\beta\mid \beta\in J\rangle$,
$W^J=\{ y\in W\mid y\beta>0\ \forall \beta\in J\}$,
and write
$xw^{-1}=y_1y_2$
with
$y_1\in W^J$, $y_2\in W_J$.
Just suppose $w\eta+\alpha$ is a weight of $E$.
Then $\nu+y_1\alpha=
y_1(\nu'+\alpha)$ would also be a
weight of $E$.
As $\nu$ is the highest weight of
$E$, 
$y_1\alpha<0$,
and hence
$\alpha\notin J$.
Then
$0<\langle\nu', \alpha^\vee\rangle=
\langle y_1\nu', y_1\alpha^\vee\rangle=\langle\nu, y_1\alpha^\vee\rangle$,
and hence
$y_1\alpha>0$, absurd.

\setcounter{equation}{0}
\noindent
(2.2)
{\bf Proposition:}
{\it
Let
$\lambda,\mu\in\Lambda$
with $\mu$ lying in the closure of the facet $\lambda$ belongs to
with respect to $W_p$.
Regarding an $L_{I,1}T$-module as a
$P_1T$-module through the quotient $P\to P/\Ru(P)$,
there is a functorial isomorphism of
$G_1T$-modules on the category of $L_{I,1}T$-modules
\[
\rT_\lambda^\mu\hat\nabla_P(?)\simeq
\hat\nabla_P(\rT_{I,\lambda}^\mu?).
\]

}

\pf
Let
$M$ be an $L_{I,1}T$-module belonging to the $\lambda$-block,
and $E$ a simple $G$-module of extremal weight
$\mu-\lambda$.
Let $\pr_\mu$ (resp. $\pr_{I,\mu}$) be the projection to the $\mu$-block of $G_1T$-
(resp. $L_{I,1}T$-)
modules.
Thus
$\rT_\lambda^\mu\hat\nabla_P(M)=\pr_\mu(E\otimes\hat\nabla_P(M))$.
If $w(\mu-\lambda)\in\Lambda_I^+$ with
$w\in W_I$
and
$v\in E\setminus0$ is of weight
$w(\mu-\lambda)$,
then
$\Dist(L_I)v$ is by (2.1)
an $L_I$-module of highest weight $w(\mu-\lambda)$.
If we put
$E'=\Dist(L_I)v$,
$
\rT_{I,\lambda}^\mu
M=\pr_{I,\mu}(E'\otimes
M)$
\cite[Remark II.7.6.1]{J}.
Thus 
\begin{align*}
\rT_\lambda^\mu\hat\nabla_P(M)
&=
\pr_\mu(E\otimes\hat\nabla_P(M))
\simeq
\pr_\mu(\hat\nabla_P(E\otimes
M))
\\
&\geq
\pr_\mu(\hat\nabla_P(E'\otimes
M))
\geq
\hat\nabla_P(\pr_{I,\mu}(E'\otimes
M))
=
\hat\nabla_P(\rT_{I,\lambda}^\mu(
M)).
\end{align*}
As it becomes an isomorphism
for $M=\hat\nabla^P(x\bullet\lambda)$ 
and
$x\in W_{I,p}$,
the isomorphism for general $M$
follows 
using the five lemma.

\setcounter{equation}{0}
\noindent
(2.3)
{\bf Corollary:}
{\it  
Let
$\lambda,\mu\in\Lambda$.
Assume that
$\mu$ 
lies in the closure of the facet $\lambda$ belongs to
with respect to
$W_p$.
Let $F_I$ be the facette $\lambda$ belongs to with respect to $W_{I,p}$ and let $\hat F_I$ be its upper closure with respect to
$W_{I,p}$.
Then
\[
\rT_\lambda^\mu
\hat\nabla_P(\hat L^P(\lambda))\simeq
\begin{cases}
\hat\nabla_P(\hat L^P(\mu))
&\text{if $\mu\in\hat F_I$},
\\
0
&\text{else},
\end{cases}
\]
in the first case of which one has a commutative diagram of
$G_1T$-modules
\[
\xymatrix{
\rT_\lambda^\mu
\hat\nabla(\lambda)
\ar[rr]^-\sim
&&
\hat\nabla(\mu)
\\
\rT_\lambda^\mu
\hat\nabla_P(\hat L^P(\lambda))
\ar@{^(->}!<0ex,2.5ex>;[u]
\ar[rr]_-\sim
&&
\hat\nabla_P(\hat L^P(\mu)).
\ar@{^(->}!<0ex,2.5ex>;[u]
}
\]

 }

\setcounter{equation}{0}
\noindent
(2.4)
For $\alpha\in R$ and $n\in\bbZ$
let
$s_{\alpha,n}$ denote the reflection in the wall
$H_{\alpha,n}$.

\begin{prop}
Let
$\lambda,\mu\in\Lambda$
with
$\lambda$ lying in an alcove
$A$ and $\mu\in\overline A$.
Assume
$
\{x\in W_p|
x\bullet\mu=\mu\}=\{e,s_{\alpha,n}\}$
for some
$\alpha\in R^+_I$ and $n\in\bbZ$.
If
$M$ is an $L_{I,1}T$-module belonging to $\mu$, there is an isomorphism of $G_1T$-modules
\[
\rT_\mu^\lambda\hat\nabla_P(M)
\simeq
\hat\nabla_P(\rT_{I,\mu}^\lambda
M),
\] 
regarding $M$ and $\rT_{I,\mu}^\lambda
M$
as
$P_1T$-modules via the quotient
$P\to
P/\Ru(P)$.

\end{prop}

\pf
Arguing as in (2.2) yields 
$\rT_\mu^\lambda\hat\nabla_P(M)
\geq
\hat\nabla_P(\rT_{I,\mu}^\lambda
M)$.
On the other hand,
if
$M=\hat\nabla^P(x\bullet\mu)$
for some
$x\in
W_{I,p}$,
\[
\ch\rT_\mu^\lambda\hat\nabla_P(\hat\nabla^P(x\bullet\mu))
=
\ch\rT_\mu^\lambda\hat\nabla(x\bullet\mu)
=
\ch\hat\nabla(x\bullet\lambda)
+\ch\hat\nabla(xs_{\alpha,n}\bullet\lambda)
\]
while
\begin{align*}
\ch\hat\nabla_P(
\rT_{I,\mu}^\lambda\hat\nabla^P(x\bullet\mu))
&=
\ch\hat\nabla_P(\hat\nabla^P(x\bullet\mu))
+\ch\hat\nabla_P(
\hat\nabla^P(xs_{\alpha,n}\bullet\lambda))
\quad\text{
as
$s_{\alpha,n}\in
W_{I,p}$}
\\
&=
\ch\hat\nabla(x\bullet\lambda)
+\ch\hat\nabla(xs_{\alpha,n}\bullet\lambda).
\end{align*}
By additivity the character equality holds for general $M$, and hence the assertion.

\setcounter{equation}{0}
\noindent
(2.5)
{\bf Corollary:}
{\it
Let
$\lambda,\mu\in\Lambda$
and keep the assumptions on
$\lambda$ and $\mu$ from (2.4).

(i)
$\rT_\mu^\lambda\hat\nabla_P(\hat L^P(\mu))$
admits a $G_1T$-filtration whose subquotients are
$\hat\nabla_P(\hat L^P(x\bullet\lambda))$,
$x\in W_{I,p}$,
with multiplicity
$m_x\in\bbN$
such that
$\ch\rT_{I,\mu}^\lambda\hat L^P(\mu)=\sum_{x\in W_{I,p}}m_x\ch\hat L^P(x\bullet\lambda)$.

(ii)
If
$\lambda<s_{\alpha,n}\bullet\lambda$, then
$\soc_{G_1T}
\rT_\mu^\lambda\hat\nabla_P(\hat L^P(\mu))=\hat
L(\lambda)$.

}

\pf
For (i) argue as in (2.2).
As
$\hat\nabla_P(\hat L^P(\mu))\leq\hat\nabla(\mu)$,
(ii) follows from the fact that
$\soc_{G_1T}\rT_\mu^\lambda\hat\nabla(\mu)=\hat L(\lambda)$.

\begin{center}
$3^\circ$
{\bf 
Grading the induction functor}
\end{center} 

In this section we employ graded representation theory from
\cite{AJS}
to show that our induction functor
$\hat\nabla_P$ can be 
graded on $p$-regular blocks.
To facilitate reference to \cite{AJS},
we will adapt to their notations 
except for
$\Bbbk=k$,
$\Lambda=X$,
and $\hat L=L_\Bbbk$.


Let
$S_\Bbbk$ be the symmetric algebra on $\bbZ R\otimes_\bbZ\Bbbk$ over
$\Bbbk$
and
$\hat S_\Bbbk$ 
its completion
along the maximal ideal
$\fm$
generated by
$R$.
We will denote
each
$\alpha\in R$ in
$S_\Bbbk$
by
$h_\alpha$ after
\cite[14.3]{AJS}.
Fix a $p$-regular weight $\lambda^+$ belonging to the bottom dominant alcove, and put
$\Omega=W_p\bullet\lambda^+$,
$Y=p\bbZ R$.
For all the unexplained notations we refer to
\cite{AJS}.

\setcounter{equation}{0}
\noindent
(3.1)
Let us first recall
\cite[\S18]{AJS}
to suit our objectives.
The category of finite
dimensional $G_1T$-modules belonging to the block
$\Omega$
may be identified with $\cC_\Bbbk(\Omega)$ from \cite{AJS}.
For each
$\lambda\in\Omega$
let
$Q_\Bbbk(\lambda)$
be the projective cover and the injective hull
of
$\hat L(\lambda)$
in
$\cC_\Bbbk(\Omega)$.
If
$Q=\bigoplus_{w\in W}Q_\Bbbk(w\bullet
\lambda^+)$,
$Q$ is a projective
$Y$-generator of
$\cC_\Bbbk(\Omega)$
\cite[E.3]{AJS}.
Thus, if
$E_{\Omega,\Bbbk}=\cC_\Bbbk(\Omega)^\sharp(Q,Q)^\op
=
\{\bigoplus_{\gamma\in Y}\cC_\Bbbk(\Omega)(Q[\gamma],Q)\}^\op$
with
$Q[\gamma]=\bigoplus_w
Q_\Bbbk(w\bullet
\lambda^++\gamma)$, it
is equipped with a structure of
finite dimensional
$Y$-graded $\Bbbk$-algebra.
Denoting the category of
$Y$-graded $E_{\Omega,\Bbbk}$-modules of finite type
by
$E_{\Omega,\Bbbk}\bfmodgr_Y$, 
the functor
$H_{\Omega,\Bbbk}=\cC_\Bbbk(\Omega)^\sharp(Q, ?)=
\bigoplus_{\gamma\in Y}\cC_\Bbbk(\Omega)(Q[\gamma], ?)$
gives an equivalence of categories from
$\cC_\Bbbk(\Omega)$ to
$E_{\Omega,\Bbbk}\bfmodgr_Y$
with quasi-inverse
$v=Q\otimes_{E_{\Omega,\Bbbk}}?$ \cite[E.4]{AJS}.

Let
$\cC(\Omega,\hat S_\Bbbk)$
denote the deformation category over
$\hat S_\Bbbk$
of
$\cC_\Bbbk(\Omega)$.
If
$\cF(\Omega,\hat S_\Bbbk)$
is its full subcategory consisting of
the objects that are free over
$\hat S_\Bbbk$,
there is
a fully faithfull functor
$\cV_\Omega$
from
$\cF(\Omega,\hat S_\Bbbk)$
to the combinatorial category
$\cK(\Omega,\hat S_\Bbbk)$
\cite[9.4]{AJS}.
Each $Q_\Bbbk(\lambda)$
lifts to a projective object
$Q_{\hat S_\Bbbk}(\lambda)$
of
$\cC(\Omega,\hat S_\Bbbk)$,
and
$\cV_\Omega
Q_{\hat S_\Bbbk}(\lambda)$
admits a graded
$S_\Bbbk$-form
$\cQ(\lambda)$
in the graded combinatorial category
$\tilde\cK(\Omega,S_\Bbbk)$.
If
$\cP=\bigoplus_{w\in W}\cQ(w\bullet
\lambda^+)$ and if
$E_\Omega=\tilde\cK(\Omega,S_\Bbbk)^\sharp(\cP,
\cP)$, then
$E_\Omega$ is a $(Y\times\bbZ)$-graded
$S_\Bbbk$-algebra of finite type
and there is an isomorphism of
$Y$-graded $\Bbbk$-algebras
$E_\Omega\otimes_{S_\Bbbk}\Bbbk\simeq
E_{\Omega,\Bbbk}$.
Thus
$E_{\Omega,\Bbbk}$ comes equipped with a structure of
finite dimensional $(Y\times\bbZ)$-graded
$\Bbbk$-algebra.
We denote by
$\tilde\cC_\Bbbk(\Omega)$ the category of
finite dimensional
$(Y\times\bbZ)$-graded
$E_{\Omega,\Bbbk}$-modules
after
\cite[18.18]{AJS}
and let
$\bar v$ denote the  functor
from
$\tilde\cC_\Bbbk(\Omega)$
to
$\cC_\Bbbk(\Omega)$
composite
of the 
forgetful functor
from
$\tilde\cC_\Bbbk(\Omega)$
to
$E_{\Omega,\Bbbk}\bfmodgr_Y$
and $v$
\cite[18.19]{AJS}.
Each $Q_\Bbbk(\lambda)$,
$Z_\Bbbk^w(\lambda\langle w\rangle)$, and
$\hat L(\lambda)$,
$\lambda\in\Omega$,
$w\in W$,
admits a graded object
$\tilde
Q_{
\Bbbk}(\lambda)$,
$\tilde
Z_{
\Bbbk}^w(\lambda)$, and
$\tilde
L_{
\Bbbk}(\lambda)$
in
$\tilde\cC_\Bbbk(\Omega
)$,
respectively,
such that 
$
\bar v\tilde Q_\Bbbk(\lambda)
\simeq
Q_\Bbbk(\lambda)$,
$
\bar v\tilde Z_\Bbbk^w(\lambda)
\simeq
Z_\Bbbk^w(\lambda\langle w\rangle)$,
and
$
\bar v\tilde L_\Bbbk(\lambda)
\simeq
L_\Bbbk(\lambda)$
in
$\cC_\Bbbk(\Omega)$
\cite[18.8 and 18.10]{AJS}.

\setcounter{equation}{0}
\noindent
(3.2)
Fix $\lambda_I^+\in
\Lambda_I^+\cap
W_p\bullet\lambda^+$
with
$\langle\lambda_I^++\rho,\alpha^\vee\rangle<p$
$\forall
\alpha\in R_I^+$.
Let
$\Omega_I=W_{I,p}\bullet\lambda_I^+$
and let $\cC_\Bbbk(\Omega_I)$ denote the category of finite dimensional
$L_{I,1}T$-modules 
belonging to the block
$\Omega_I$.
Put
$Y_I=p\bbZ I$.
For each
$\nu\in\Lambda$ let
$Q_{I,\Bbbk}(\nu)$ be the projective cover of $\hat L^P(\nu)$ as $L_{I,1}T$-module.
If
$Q_I=\bigoplus_{w\in W_I}Q_{I,\Bbbk}(w\bullet\lambda_I^+)$,
it is a projective $Y_I$-generator of $\cC_\Bbbk(\Omega_I)$.
Let
$E_{\Omega_I,\Bbbk}
=\cC_\Bbbk(\Omega_I)^\sharp(Q_{I,\Bbbk}, Q_{I,\Bbbk})^\op$.
It is equipped with a structure of finite
dimensional
$(Y_I\times\bbZ)$-graded $\Bbbk$-algebra.
Let
$\tilde\cC_\Bbbk(\Omega_I)$
denote the category of
$(Y_I\times\bbZ)$-graded
$E_{\Omega_I,\Bbbk}$-modules,
and construct
$\tilde\nabla_{I,\Bbbk}(\lambda),\tilde
L_{I,\Bbbk}(\lambda)\in
\tilde\cC_\Bbbk(\Omega_I)$,
$\lambda\in\Omega_I$,
just like 
$\tilde\nabla_{\Bbbk}(\lambda)$,
$\tilde
L_{\Bbbk}(\lambda)$
for $G$.

Unless otherwise specified we will regard an $L_{I,1}T$-module as a $P_1T$-module via 
inflation along the quotient
$P\to
P/\Ru(P)\simeq
L_I$.

\begin{lem}
There is a functorial isomorphism
from the category of 
$Y_I$-graded
$E_{\Omega_I,\Bbbk}$-modules
of finite type to
$\cC_\Bbbk(\Omega)$
\[
Q\otimes_{E_{\Omega,\Bbbk}}\cC_\Bbbk(\Omega)^\sharp(Q,\hat\nabla_P(Q_I))\otimes_{E_{\Omega_I,\Bbbk}}?
\simeq
\hat\nabla_P(Q_I\otimes_{E_{\Omega_I,\Bbbk}}?).\]

\end{lem}

\pf
Let
$\tilde M$ be a
$Y_I$-graded
$E_{\Omega_I,\Bbbk}$-module
of finite type.
As
$Q\otimes_{E_{\Omega,\Bbbk}}?$
and
$\cC_\Bbbk(\Omega)^\sharp(Q, ?)$
are quasi-inverse to each other,
$Q\otimes_{E_{\Omega,\Bbbk}}\cC_\Bbbk(\Omega)^\sharp(Q,\hat\nabla_P(Q_I))\otimes_{E_{\Omega_I,\Bbbk}}\tilde M\simeq
\hat\nabla_P(Q_I)\otimes_{E_{\Omega_I,\Bbbk}}\tilde M$,
which is isomorphic to
$\hat\nabla_P(Q_I\otimes_{E_{\Omega_I,\Bbbk}}\tilde M)$
if
$\tilde M$
is isomorphic to
$E_{\Omega_I,\Bbbk}$.
In general, apply the five lemma
to a natural homomorphism of
$G_1T$-modules
$\hat\nabla_P(Q_I)\otimes_{E_{\Omega_I,\Bbbk}}\tilde M
\to
\hat\nabla_P(Q_I\otimes_{E_{\Omega_I,\Bbbk}}\tilde M)$.

\setcounter{equation}{0}
\noindent
(3.3)
We will show that
the lemma above refines to a
commutative diagram
\begin{equation}
\xymatrix{
\tilde\cC_\Bbbk(\Omega_I)
\ar[rrrr]^-{\cC_\Bbbk(\Omega)^\sharp(Q,\hat\nabla_P(Q_I))\otimes_{E_{\Omega_I,\Bbbk}}?}
\ar[d]_-{\bar v_I}
&&&&
\tilde\cC_\Bbbk(\Omega)
\ar[dd]^-{\bar v}
\\
L_{I,1}T\Mod
\ar[d]
&&&&
\\
P_{1}T\Mod
\ar[rrrr]_{\hat\nabla_P}&&&&
G_{1}T\Mod
}
\end{equation}
in such a way that for each
$\lambda\in\Omega_I$
\begin{equation}
\cC_\Bbbk(\Omega)^\sharp(Q,
\hat\nabla_P(Q_I))\otimes_{E_{\Omega_I,\Bbbk}}
\tilde\nabla_{I,\Bbbk}(\lambda)\simeq
\tilde\nabla_\Bbbk(\lambda)\langle\delta(\lambda)-\delta_I(\lambda)\rangle,
\end{equation}
where
$\delta$
(resp.
$\delta_I$)
is the length function on
$\Omega$
(resp.
$\Omega_I$)
\cite[17.1]{AJS}.

To justify the commutative diagram, we have only to show that
$\cC_\Bbbk(\Omega)^\sharp(Q,
\hat\nabla_P(Q_I))$ is equipped with a structure of
$(Y\times\bbZ)$-graded left
$E_{\Omega,\Bbbk}$
and $(Y_I\times\bbZ)$-graded right $E_{\Omega_I,\Bbbk}$-bimodule.
For that we first
deform the functor
$\hat\nabla_P$.
Put $S_{I,\Bbbk}=S_\Bbbk(\bbZ 
I\otimes_\bbZ\Bbbk)$
to be
the symmetric algebra
over
$\Bbbk$ on
$\bbZ 
I\otimes_\bbZ\Bbbk$.
We will write
$A_G$
(resp. $A_I$) 
for $\hat S_\Bbbk$ (resp.
the completion of $S_{I,\Bbbk}$ with respect to the maximal ideal generated by $\bbZ 
I\otimes_\bbZ\Bbbk$).
For each $\beta\in R^+$
let
$A_G^\beta=A_G[\frac{1}{h_\alpha}\mid\alpha\in R^+\setminus\{\beta\}]$,
$A_G^\emptyset=A_G[\frac{1}{h_\alpha}\mid\alpha\in R^+]$,
and for
$\beta\in R_I^+$
put
$A_I^{\beta}=A_I[\frac{1}{h_\alpha}\mid\alpha\in R_I^+\setminus\{\beta\}]$,
$A_I^{\emptyset}=A_I[\frac{1}{h_\alpha}\mid\alpha\in R_I^+]$.
We will regard $A_G$ as an $A_I$-algebra via 
inclusion $R_I\hookrightarrow
R$;
in case $R$ has two lengths, if a component $I'$ of $I$ consists only of long roots, we 
take
$h_\alpha=d_\alpha
H_\alpha$
for each $\alpha\in R_{I'}$ with
$d_\alpha$ for $R$
instead of $h_\alpha=H_\alpha$.
Though this deviates from the convention in
\cite[14.4/p. 11]{AJS}, it causes no difference to our application.
Thus $A_G^\emptyset$ is an $A_I^\emptyset$-algebra,
and for
$\beta\in R^+$
\[
A_G^\beta
\simeq
\begin{cases}
A_I^\beta\otimes_{A_I
}A_G^\beta
&\text{
if
$\beta\in R_I^+$}
\\
A_I^\emptyset\otimes_{A_I
}A_G^\beta
&\text{
else}.
\end{cases}
\]

For a $W_{I,p}$-orbit
$\Gamma_I$ in $\Lambda$
define
$\cC(\Gamma_I, A_I)$,
$\cC(\Gamma_I, A_I^\beta)$
for
$\beta\in R_I^+$,
$\cC(\Gamma_I, A_G)$,
$\cC(\Gamma_I, A_G^\beta)$
for
$\beta\in R^+$,
$\cF\cC(\Gamma_I, A_I)$,
$\cK(\Gamma_I, A_I)$,
$\tilde\cK(\Gamma_I, S_{I,\Bbbk})$ and $\tilde\cC_\Bbbk(\Gamma_I)$
for $L_I$ just as for $G$;
precisely these are defined first
for the semisimple part of $L_I$ 
and then extended to $L_I$
in a natural way.
For $\nu\in\Lambda$
define likewise
$Z_{I,A_I}(\nu)$,
$Z_{I,A_I}^\beta(\nu)
=
Z_{I,A_I^{\beta}}(\nu)$ for
$\beta\in R_I^+$,
and
$Z_{I,A_I}^\emptyset(\nu)
=
Z_{I,A_I^{\emptyset}}(\nu)$
as well as
$Z_{I,A_G}(\nu)$,
$Z_{I,A_G^{\beta}}(\nu)$ for
$\beta\in R^+$,
and
$Z_{I,A_G^{\emptyset}}(\nu)$.

Recall from (1.7) the parabolic subgroup
$P^+=\langle B^+, U_{-\alpha}|\alpha\in I\rangle$.
Let
$\Gamma$ be the
$W_p$-orbit in $\Lambda$ containing
$\Gamma_I$.
Regarding an object of
$\cF\cC(\Gamma_I, A_I)$
as a 
$\Dist(P^+_1)$-module
by the quotient
$P^+\to
P^+/\Ru(P^+)$,
define a functor 
$\hat\nabla_{P,A_I}: \cF\cC(\Gamma_I, A_I)\to
\cF\cC(\Gamma, A_G)$ via
\[
M\mapsto
(\Dist(G_1)\otimes_{\Dist(P^+_1)}M^\tau)^\tau
\otimes_{A_I}A_G\simeq
\{
\Dist(G_1)\otimes_{\Dist(P^+_1)}(M\otimes_{A_I}A_G)^\tau\}^\tau,
\]
which reduces to 
$\hat\nabla_P$ by reduction to $\Bbbk$.
For each $\nu\in\Gamma_I$
one has
\begin{equation}
\hat\nabla_{P,A_I}(Z_{I,A_I}(\nu)^\tau)
\simeq
Z_{A_G}(\nu)^\tau.
\end{equation}

\setcounter{equation}{0}
\noindent
(3.4)
Let $U_1(w_I)=\prod_{\beta\in R^+\setminus R_I}U_{-\beta,1}$
be the Frobenius kernel of the unipotent radical of
$P$
and
$\Dist^+(U_1(w_I))$ the augmentation ideal of
$\Dist(U_1(w_I))$.
Let
$\Gamma$ be an arbitrary
$W_p$-orbit.
For each 
$M\in\cC(\Gamma, A_G)$ put
\[
M_\frn=M/\Dist^+(U_1(w_I))M
\simeq
\{
\Dist(U_1(w_I))/\Dist^+(U_1(w_I))\}\otimes_{\Dist(U_1(w_I))}
M
\]
the module of $\Dist^+(U_1(w_I))$-coinvariants of $M$.
If $M=Z_{A_G}(\nu)$, $\nu\in\Lambda$,
taking the $\tau$-dual of (1.7) yields
an isomorphism in
$\cC(W_{I,p}\bullet\nu, A_G)$
\begin{equation}
Z_{A_G}(\nu)_\frn\simeq
Z_{I,A_I}(\nu)\otimes_{A_I}A_G
\simeq
Z_{I,A_G}(\nu).
\end{equation}

Let
$\beta\in R_I^+$, $\nu\in\Gamma$ with
$\beta\uparrow\nu>\nu$,
and put
$\Gamma_I=W_{I,p}\bullet\nu$.
One has from
\cite[8.6]{AJS},
as
$d_\beta\in\Bbbk^\times$ by the standing hypothesis on
$p$
\cite[14.4]{AJS},
\begin{align}
\Ext^1_{\cC(\Gamma,A_G^\beta)}
(Z_{A_G}^\beta(\nu), 
&
Z_{A_G}^\beta(\beta\uparrow\nu))
\simeq
A_G^\beta
h_\beta^{-1}/A_G^\beta
\simeq
(A_I^\beta
h_\beta^{-1}/A_I^\beta)\otimes_{A_I}A_G^\beta
\\
\notag
&\simeq
\Ext^1_{\cC(\Gamma_I,A_I^\beta)}
(Z_{I,A_I}^\beta(\nu), 
Z_{I,A_I}^\beta(\beta\uparrow\nu))
\otimes_{A_I}A_G^\beta
\\
\notag
&\simeq
\Ext^1_{\cC(\Gamma_I,A_G^\beta)}
(Z_{I,A_G^\beta}(\nu), 
Z_{I,A_G^\beta}(\beta\uparrow\nu))
\quad\text{by \cite[3.2]{AJS}}.
\end{align}

\begin{lem}
Assume
$\beta\uparrow\nu>\nu$.
If
$0\to
Z_{A_G}^\beta(\beta\uparrow\nu)\to
M\to
Z_{A_G}^\beta(\nu)\to0$
is exact in
$\cC(\Gamma, A_G^\beta)$,
applying
$?_\frn$ to the sequence
yields
an exact sequence
$0\to
Z_{I,A_G^\beta}(\beta\uparrow\nu)\to M_\frn\to
Z_{I,A_G^\beta}(\nu)\to0$ 
with $M$ projective in $\cC(\Gamma, A_G^\beta)$
iff
$M_\frn$ projective in
$\cC(\Gamma_I, A_G^\beta)$.
Conversely,
applying $\Dist(G_1)\otimes_{\Dist(
P_1^+)}?$
to the latter sequence recovers the former.
Likewise,
if
$0\to
Z_{I,A_I}^\beta(\beta\uparrow\nu)\to
M'\to
Z_{I,A_I}^\beta(\nu)\to0$
is an exact sequence
in
$\cC(\Gamma_I, A_I^{\beta})$
with
$M'$ projective,
then applying 
$\Dist(G_1)\otimes_{\Dist(
P_1^+)}?\otimes_{A_I}A_G^\beta$
yields an exact sequence
$0\to
Z_{A_G}^\beta(\beta\uparrow\nu)\to
\Dist(G_1)\otimes_{\Dist(
P_1^+)}M'
\otimes_{A_I}A_G^\beta
\to
Z_{A_G}^\beta(\nu)\to0$
with
$\Dist(G_1)\otimes_{\Dist(
P_1^+)}M'
\otimes_{A_I}A_G^\beta$
projective
in
$\cC(\Gamma, A_G^\beta)$.

\end{lem}

\pf
Assume the sequence
$0\to
Z_{A_G}^\beta(\beta\uparrow\nu)\to
M\to
Z_{A_G}^\beta(\nu)\to0$
is exact.
As
$?_\frn\simeq
\{
\Dist(U_1(w_I))/\Dist^+(U_1(w_I))\}\otimes_{\Dist(U_1(w_I))}
?$
and as
$Z_{A_G}^\beta(\nu)\simeq
\Dist(U_1)\simeq
\Dist(U_1(w_I))\otimes_\Bbbk
\Dist((B\cap L_I)_1)$
is free over
$\Dist(U_1(w_I))$,
$0\to
Z_{A_G}^\beta(\beta\uparrow\nu)_\frn\to M_\frn\to
Z_{A_G}^\beta(\nu)_\frn\to0$
remains exact
with
$Z_{A_G}^\beta(\beta\uparrow\nu)_\frn\simeq
Z_{I,A_G^\beta}(\beta\uparrow\nu)$
and
$Z_{A_G}^\beta(\nu)_\frn\simeq
Z_{I,A_G^\beta}(\nu)$.

Recall from
\cite[12.4]{AJS} how each $M$ is constructed.
Let $w_\beta\in W_I$ with
$w_\beta^{-1}\beta\in I$.
Let
$v_\nu^{w_\beta}\in
Z_{A_G}^\beta(\nu)$ of weight
$\nu\langle w_\beta\rangle$
corresponding to
the standard generator $1\otimes1$
of $Z_{A_G^\beta}^{w_\beta}(\nu\langle w_\beta\rangle)$
under the isomorphism
$Z_{A_G}^\beta(\nu)=
Z_{A_G^\beta}(\nu)\simeq
Z_{A_G^\beta}^{w_\beta}(\nu\langle w_\beta\rangle)$,
and define
$v_{\beta\uparrow\nu}^{w_\beta}\in
Z_{A_G}^\beta(\beta\uparrow\nu)$ likewise.
Write
$\langle\nu+\rho, \beta^\vee\rangle\equiv p-n\mod p$
with
$n\in[0, p[$,
and put
$z_\nu=
E_{-\beta}^{(n)}v_{\beta\uparrow\nu}^{w_\beta}b+v_\nu^{w_\beta}\in
Z_K(\beta\uparrow\nu)\oplus
Z_K(\nu)$ for each
$b\in A_G^\beta h_\beta^{-1}$
with
$K=\Frac(A_G)$, so
$z_\nu$ is of weight
$\nu\langle w_\beta\rangle$.
Then $M$ is of the form
$M_\nu^{w_\beta}(b)=
\Dist(G_1)v_{\beta\uparrow\nu}^{w_\beta}A_G^\beta+
\Dist(G_1)z_\nu
A_G^\beta$
living in
$Z_K(\beta\uparrow\nu)\oplus
Z_K(\nu)$, and the sequence reads
$v_{\beta\uparrow\nu}^{w_\beta}$ mapping to itself while
$z_\nu\mapsto
v_\nu^{w_\beta}$.
Now
\begin{align*}
\Dist(G_1)
&\simeq
\Dist(^{w_\beta}U_1)\otimes
\Dist(^{w_\beta}B_1^+)
\\
&\simeq
\Dist(^{w_\beta}U_1(w_I))\otimes
\Dist(^{w_\beta}(B\cap
L_I)_1)
\otimes
\Dist(^{w_\beta}B_1^+)
\\
&\simeq
\Dist(U_1(w_I))\otimes
\Dist(^{w_\beta}(B\cap
L_I)_1)
\otimes
\Dist(^{w_\beta}B_1^+)
\end{align*}
as
$^{w_\beta}U(w_I)=^{w_\beta}\prod_{\alpha\in R^+\setminus R_I}U_{-\alpha}=
\prod_{\alpha\in R^+\setminus R_I}U_{-\alpha}=
U(w_I)$.
Thus
$(\Dist(G_1)v_{\beta\uparrow\nu}^{w_\beta})_\frn
\simeq
\Dist(^{w_\beta}(B\cap
L_I)_1)v_{\beta\uparrow\nu}^{w_\beta}$,
$(\Dist(G_1)z_\nu)_\frn\simeq
\Dist(^{w_\beta}(B\cap
L_I)_1)z_\nu$,
and hence
$M_\nu^{w_\beta}(b)_\frn=
\Dist(L_{I,1})v_{\beta\uparrow\nu}^{w_\beta}A_G^\beta+
\Dist(L_{I,1})z_\nu
A_G^\beta$.
It follows from \cite[8.7]{AJS}
that
$M_\nu^{w_\beta}(b)$ is projective in
$\cC(\Omega, A_G^\beta)$
iff $A_G^\beta b=A_G^\beta
h_\beta^{-1}
/A_G^\beta$ iff
$M_\nu^{w_\beta}(b)_\frn$
is projective
in
$\cC(\Omega_I, A_G^\beta)$.
Likewise the last assertion follows from (1).

\setcounter{equation}{0}
\noindent
(3.5)
We now transfer from $\cF\cC(\Omega,A_G)$
(resp.
$\cF\cC(\Omega_I,A_I)$)
to the combinatorial category
$\cK(\Omega,A_G)$
(resp.
$\cK(\Omega_I,A_I)$)
via the fully faithful functor
$\cV_\Omega$
(resp.
$\cV_{\Omega_I}$).
Define a functor $\cI:
\cK(\Omega_I, A_I)\to\cK(\Omega, A_G)$ as follows: for each
$\cM\in\cK(\Omega_I, A_I)$ 
and $\lambda\in\Omega$
set 
\[
(\cI\cM)(\lambda)
=
\begin{cases}
\cM(\lambda)\otimes_{A_I
}
A_G^\emptyset
&\text{if
$\lambda\in\Omega_I$}
\\
0
&\text{else},
\end{cases}
\]
and for each $\beta\in R^+$
set
\[
(\cI\cM)(\lambda,\beta)
=
\begin{cases}
\cM(\lambda,\beta)\otimes_{A_I
}
A_G^\beta
&\text{if
$\lambda\in\Omega_I$
and $\beta\in R_I^+$}
\\
\cM(\lambda)\otimes_{A_I
}
A_G^\beta
&\text{if
$\lambda\in\Omega_I$
and $\beta\notin R_I^+$}
\\
\cM(\beta\uparrow\lambda)\otimes_{A_I
}
A_G^\beta
&\text{if
$\beta\uparrow\lambda\in\Omega_I$
and $\beta\notin R_I^+$}
\\
0
&\text{else}.
\end{cases}
\]
Let
$\cD(\Omega, A_G)$
(resp.
$\cD(\Omega_I, A_I)$)
be the full subcategory of
$\cC(\Omega, A_G)$
(resp.
$\cC(\Omega_I, A_I)$)
consisting of objects admitting a
$Z_{A_G}$-
(resp.
$Z_{I,A_I}$-)
filtration.
We want to show a functorial isomorphism
$\cV_\Omega\circ
\hat\nabla_{P,A_I}\simeq
\cI\circ\cV_{\Omega_I}$
from
$\cD(\Omega_I, A_I)$
to
$\cD(\Omega, A_G)$.

Recall that
$\cV_\Omega$ and 
$\cV_{\Omega_I}$ are defined with specific choice of extensions
according to
Theorem of Good Choices
\cite[13.4]{AJS}.

\begin{lem}
Let
$\lambda\in\Omega_I$ and $\beta\in R_I^+$ with
$\beta\uparrow\lambda>\lambda$.
Let
$e^\beta(\lambda)\in\Ext_{\cC(\Omega,A_G^\beta)}^1(Z_{A_G}^\beta(\lambda),
Z_{A_G}^\beta(\beta\uparrow\lambda))$
and
$e_I^\beta(\lambda)\in\Ext_{\cC(\Omega_I,A_I^\beta)}^1(Z_{I,A_I}^\beta(\lambda),
Z_{I,A_I}^\beta(\beta\uparrow\lambda))$
chosen according to Theorem of Good Choices.
Let
$Y_{A_G}^\beta(\lambda)\in
\cC(\Omega,A_G^\beta)$
(resp. 
$Y_{I,A_I}^\beta(\lambda)\in
\cC(\Omega_I,A_I^\beta)$)
be the module representing
$e^\beta(\lambda)$
(resp. 
$e_I^\beta(\lambda)$).
Then
$Y_{I,A_I}^\beta(\lambda)\otimes_{A_I}A_G^\beta=
Y_{A_G}^\beta(\lambda)_\frn$
and
$\Dist(G_1)\otimes_{\Dist(P^+_1)}Y_{I,A_I}^\beta(\lambda)\otimes_{A_I}A_G^\beta=
Y_{A_G}^\beta(\lambda)$.

\end{lem}

\pf
Write
$\lambda=w_1\bullet\lambda^++p\gamma_1=
w_2\bullet\lambda_I^++p\gamma_2$
with
$w_1\in W$, $w_2\in W_I$,
$\gamma_1\in\bbZ R$,
$\gamma_2\in\bbZ R_I$.
By \cite[13.25]{AJS} we may assume $\lambda_I^+=w_2^{-1}w_1\bullet\lambda^+$.
Then for each $\alpha\in R_I^+$
\begin{equation}
(w_2^{-1}w_1)^{-1}\alpha=w_1^{-1}w_2\alpha>0,
\end{equation}
and hence
\begin{equation}
w_1^{-1}\alpha>0
\quad\text{ iff }\quad
w_2^{-1}\alpha>0.
\end{equation}

Recall from \cite[13.2.5]{AJS}
that
$e^\beta(\lambda)=b^\beta(\lambda)e_0^\beta(\lambda)$,
and thus 
$e_I^\beta(\lambda)=b_I^\beta(\lambda)e_{I,0}^\beta(\lambda)$
likewise
with the RHS's specified as follows.
By \cite[13.2.3, 4]{AJS}
\begin{equation}
b^\beta(\lambda)
=
\begin{cases}
\varepsilon^\beta_{w_1\bullet\lambda^+,-\rho}
d(w_1\bullet\lambda^+,-\rho,s_\beta)\kappa(\beta)
&
\text{if $w_1^{-1}\beta>0$}
\\
\varepsilon^\beta_{w_1\bullet\lambda^+,-\rho}
d(w_1\bullet\lambda^+,-\rho,s_\beta)h_\beta
&
\text{else}
\end{cases}
\end{equation}
with
$\kappa(\beta)=\displaystyle
\underset{\substack{
\alpha\in R^+
\\
s_\beta\alpha<0,
w_\beta^{-1}\alpha<0}}{\prod}
h_\alpha^{-\langle\beta,\alpha^\vee\rangle}$,
where
$w_\beta\in W_I$ such that $w_\beta^{-1}\beta\in I$, and thus
\begin{equation}
b_I^\beta(\lambda)
=
\begin{cases}
\varepsilon^\beta_{I,w_2\bullet\lambda_I^+,-\rho_I}
d_I(w_2\bullet\lambda_I^+,-\rho_I,s_\beta)\kappa_I(\beta)
&
\text{if $w_2^{-1}\beta>0$}
\\
\varepsilon^\beta_{I,w_2\bullet\lambda_I^+,-\rho_I}
d_I(w_2\bullet\lambda_I^+,-\rho_I,s_\beta)h_\beta
&
\text{else}
\end{cases}
\end{equation}
with
$\kappa_I(\beta)=\displaystyle
\underset{\substack{
\alpha\in R_I^+
\\
s_\beta\alpha<0,
w_\beta^{-1}\alpha<0}}{\prod}
h_\alpha^{-\langle\beta,\alpha^\vee\rangle}$.
By (2) the two cases in (3) and (4) agree,
and
$\kappa(\beta)=\kappa_I(\beta)$.
By \cite[12.12.5]{AJS}
\[
\varepsilon^\beta_{w_1\bullet\lambda^+,-\rho}
=
\underset{\substack{
\alpha\in R^+
\\
s_\beta\alpha<0}}{\prod}
(-1)^{\langle-\rho-w_1\bullet\lambda^+,\alpha^\vee\rangle
\bar\alpha(-\rho-w_1\bullet\lambda^+)}
\underset{\substack{
\alpha\in R^+\setminus\{\beta\}
\\
s_\beta\alpha<0,w_\beta^{-1}\alpha>0}}{\prod}
(-1)^{\langle-\rho-w_1\bullet\lambda^+,\alpha^\vee\rangle}
\]
with
$\bar\alpha(\nu)=
\begin{cases}
1
&\text{if
$\langle\nu,\alpha^\vee\rangle>0$}
\\
0
&\text{else}
\end{cases}
$
for each
$\nu\in\Lambda$
\cite[A.1.1]{AJS}, and thus
\[
\varepsilon^\beta_{I,w_2\bullet\lambda_I^+,-\rho_I}
=
\underset{\substack{
\alpha\in R_I^+
\\
s_\beta\alpha<0}}{\prod}
(-1)^{\langle-\rho_I-w_2\bullet\lambda_I^+,\alpha^\vee\rangle
\bar\alpha(-\rho_I-w_2\bullet\lambda_I^+)}
\underset{\substack{
\alpha\in R_I^+\setminus\{\beta\}
\\
s_\beta\alpha<0,w_\beta^{-1}\alpha>0}}{\prod}
(-1)^{\langle-\rho_I-w_2\bullet\lambda_I^+,\alpha^\vee\rangle}.
\]
One has
$-\rho-w_1\bullet\lambda^+
=
-\rho-w_1\bullet
(w_2^{-1}w_1)^{-1}\bullet\lambda_I^+
=-\rho_P-\rho_I-w_2\bullet\lambda_I^+$.
If $\alpha\in R^+$ with $s_\beta\alpha<0$, $\alpha\in R_I^+$, and hence
$\langle
-\rho-w_1\bullet\lambda^+,\alpha^\vee\rangle
=
\langle
-\rho_P-\rho_I-w_2\bullet\lambda_I^+,\alpha^\vee\rangle
=
\langle
-\rho_I-w_2\bullet\lambda_I^+,\alpha^\vee\rangle$
and
$\bar\alpha(
-\rho-w_1\bullet\lambda^+)
=
\bar\alpha(
-\rho_I-w_2\bullet\lambda_I^+)$.
Thus
$\varepsilon^\beta_{w_1\bullet\lambda^+,-\rho}=
\varepsilon^\beta_{I,w_2\bullet\lambda_I^+,-\rho_I}$.
By
\cite[13.2.2]{AJS}
\begin{align*}
d(w_1\bullet\lambda^+,-\rho,s_\beta)
&=
\underset{\substack{
\alpha\in R^+
\\
s_\beta\alpha<0,w_1^{-1}\alpha<0}}{\prod}
\frac{[k_\alpha;w_1\bullet\lambda^++\rho]}{h_\alpha},
\\
d_I(w_2\bullet\lambda_I^+,-\rho_I,s_\beta)
&=
\underset{\substack{
\alpha\in R_I^+
\\
s_\beta\alpha<0,w_2^{-1}\alpha<0}}{\prod}
\frac{[k_\alpha;w_2\bullet\lambda_I^++\rho_I]_I}{h_\alpha}.
\end{align*}
By (2) again the products run over the same subset of $R_I^+$.
By
\cite[13.1.4]{AJS}
\begin{align*}
[k_\alpha;w_1\bullet\lambda^++\rho]
&=
(H_\alpha+\langle
w_1\bullet\lambda^++\rho,\alpha^\vee\rangle)H_\alpha^{-1}h_\alpha
\end{align*}
with
$\langle
w_1\bullet\lambda^++\rho,\alpha^\vee\rangle
=
\langle
w_1\bullet(w_2^{-1}w_1)^{-1}\bullet\lambda_I^++\rho,\alpha^\vee\rangle
=\langle w_2\bullet\lambda_I^++\rho,\alpha^\vee\rangle
=
\langle
w_2\bullet\lambda_I^++\rho_I,\alpha^\vee\rangle$,
and hence
$d(w_1\bullet\lambda^+,-\rho,s_\beta)
=
d_I(w_2\bullet\lambda_I^+,-\rho_I,s_\beta)$ 
and
$b^\beta(\lambda)=b_I^\beta(\lambda)$.

We compare 
next
$e_0^\beta(\lambda)$
and
$e_{I,0}^\beta(\lambda)$.
Take
$\omega\in\Lambda$ in the upper closure of the facet
$\lambda$ belongs to with respect to
$\langle
s_{\beta,r}
\mid
r\in\bbZ\rangle$.
By
\cite[12.13.1]{AJS}
\begin{align}
t_0^\beta[\omega,\lambda]
e_0^\beta(\lambda)
&=
\varepsilon^\beta_{\lambda\omega}
d(\omega,\lambda,s_\beta)h_\beta^{-1}+A_G^\beta,
\\
t_{I,0}^\beta[\omega,\lambda]
e_{I,0}^\beta(\lambda)
\notag&=
\varepsilon^\beta_{I,\lambda,\omega}
d_I(\omega,\lambda,s_\beta)h_\beta^{-1}+A_I^\beta.
\end{align}
By definition
\cite[12.12.5]{AJS}
again
\begin{align}
\varepsilon_{\lambda\omega}^\beta
&=
\underset{\substack{
\alpha\in R^+
\\
s_\beta\alpha<0}}{\prod}
(-1)^{\langle\omega-\lambda,\alpha^\vee\rangle
\bar\alpha(\omega-\lambda)}
\underset{\substack{
\alpha\in R^+
\\
s_\beta\alpha<0,
w_\beta^{-1}\alpha>0}}{\prod}
(-1)^{\langle\omega-\lambda,\alpha^\vee\rangle
}
=
\varepsilon_{I,\lambda,\omega}^\beta.
\end{align}
By definition
\cite[A.7.1 and A.2.1]{AJS}
\[
d(\omega,\lambda,s_\beta)
=
\underset{\substack{
\alpha\in R^+
\\
s_\beta\alpha<0}}{\prod}
d(\omega,\lambda,\alpha)
=
\underset{\substack{
\alpha\in R^+
\\
s_\beta\alpha<0}}{\prod}
(\frac{H_\alpha+\langle\omega+\rho,\alpha^\vee\rangle}
{H_\alpha+\langle\lambda+\rho,\alpha^\vee\rangle}
)^{\bar\alpha(\omega-\rho)}
=
d_I(\omega,\lambda,s_\beta).
\]
It follows in (5)
that
$\varepsilon^\beta_{\lambda\omega}
d(\omega,\lambda,s_\beta)h_\beta^{-1}=
\varepsilon^\beta_{I,\lambda,\omega}
d_I(\omega,\lambda,s_\beta)h_\beta^{-1}$.
By \cite[12.12.1]{AJS}
\[
t_0^\beta[\omega,\lambda]=
t[\omega,\lambda,a_{\lambda\omega}],
\quad
t_{I,0}^\beta[\omega,\lambda]=
t_I[\omega,\lambda,a_{I,\lambda,\omega}]
\]
with
\begin{align}
a_{\lambda\omega}
&=
a_{\lambda\omega}'
\varepsilon_{\lambda\omega}^\beta
\\
\notag&=
\varepsilon_{\lambda\omega}^\beta
\quad\text{by \cite[A.12]{AJS}}
\\
\notag&=
\varepsilon_{I,\lambda,\omega}^\beta
\quad\text{by (6)}
\\
\notag&=
a_{I,\lambda,\omega}'
\varepsilon_{I,\lambda,\omega}^\beta
=
a_{I,\lambda,\omega}.
\end{align}
We have
$t[\omega,\lambda,a_{\lambda\omega}]
=
t[\omega,\lambda,e,\bar e]$
by \cite[12.8.2]{AJS}
with
$e\in E_{\omega-\lambda}\setminus0$,
$E$ a simple $G$-module of extremal weight $\omega-\lambda$
\cite[11.1]{AJS}, and with
$\bar e=a_{\lambda\omega}(-1)^nE_{-\beta}^{(n)}e\in
E_{s_\beta(\omega-\lambda)}\setminus0$,
$\langle\lambda+\rho,\beta^\vee\rangle\equiv p-n
\mod p$,
$n\in]0,p[$,
by
\cite[12.3.1]{AJS}.
Recall from
\cite[12.6]{AJS}
the definition of
$t[\omega,\lambda,e,\bar e]:
\Ext^1_{\cC(\Omega,A_G^\beta)}(Z_{A_G}^\beta(\lambda),
Z_{A_G}^\beta(\beta\uparrow\lambda))
\to
H_\beta^{-1}A_G^\beta/A_G^\beta=
h_\beta^{-1}A_G^\beta/A_G^\beta$.
Let
$\xi\in
\Ext^1_{\cC(\Omega,A_G^\beta)}(Z_{A_G}^\beta(\lambda),
Z_{A_G}^\beta(\beta\uparrow\lambda))$
represented by
a short exact sequence
\begin{equation}
0\to
Z_{A_G}^\beta(\beta\uparrow\lambda)
\overset{i}{\to}
M\overset{j}{\to}
Z_{A_G}^\beta(\lambda)
\to
0.
\end{equation}
As $H_\beta\xi=0$,
there is
$j'\in\cC(\Omega,A_G^\beta)(Z_{A_G}^\beta(\lambda),M)$
with
$j\circ
j'=H_\beta\id_{Z_{A_G}^\beta(\lambda)}$.
Apply the translation functor
$T_\lambda^\omega$ to (8) to obtain
a split exact sequence
\[
0\to
T_\lambda^\omega
Z_{A_G}^\beta(\beta\uparrow\lambda)
\overset{T_\lambda^\omega
i}{\longrightarrow}
T_\lambda^\omega
M\overset{T_\lambda^\omega
j}{\longrightarrow}
T_\lambda^\omega
Z_{A_G}^\beta(\lambda)
\to
0.
\]
Let
$i'\in\cC(A_G^\beta)(T_\lambda^\omega
Z_{A_G}^\beta(\lambda),T_\lambda^\omega
M)$
with
$i'\circ
T_\lambda^\omega
j=
\id_{T_\lambda^\omega
Z_{A_G}^\beta(\lambda)}$.
Recall from
\cite[11.2.1]{AJS}
isomorphisms
$f_e: Z_{A_G}^\beta(\omega)\to
T_\lambda^\omega
Z_{A_G}^\beta(\lambda)
=
\pr(E\otimes
Z_{A_G}^\beta(\lambda))$
via
$1\otimes1\mapsto
\pr(e\otimes1\otimes1)$
and
$f_{\bar e}: Z_{A_G}^\beta(\omega)\to
T_\lambda^\omega
Z_{A_G}^\beta(\beta\uparrow
\lambda)
=
\pr(E\otimes
Z_{A_G}^\beta(\beta\uparrow\lambda))$
via
$1\otimes1\mapsto
\pr(\bar e\otimes1\otimes1)$.
If $a\in
A_G^\beta$ with
$f_{\bar e}^{-1}\circ
i'\circ
T_\lambda^\omega
j'\circ f_e=a\id_{Z_{A_G}^\beta(\omega)}$, then
$t[\omega,\lambda,e,\bar e]\xi=aH_\beta^{-1}+A_G^\beta$.
Now recall  
the $L_I$-submodule $E'$
of $E$ 
from (2.2)
and choose
$e_I=e\in E'$ and
$\bar e_I=\bar e\in E'$
to define
$t_I[\omega,\lambda,e_I,\bar e_I]:
\Ext^1_{\cC(\Omega_I,A_I^\beta)}(Z_{I,A_I}^\beta(\lambda),
Z_{I,A_I}^\beta(\beta\uparrow\lambda))
\to
H_\beta^{-1}A_I^\beta/A_I^\beta=
h_\beta^{-1}A_I^\beta/A_I^\beta$ likewise.
As we have natural isomorphisms
from (2.2) or rather from its $\tau$-dual
\begin{align*}
\Dist(G_1)\otimes_{\Dist(P_1^+)}Z_{I,A_I}^\beta(\lambda)\otimes_{A_I}A_G^\beta
&\simeq
Z_{A_G}^\beta(\lambda),
\\
\Dist(G_1)\otimes_{\Dist(P_1^+)}Z_{I,A_I}^\beta(\beta\uparrow\lambda)\otimes_{A_I}A_G^\beta
&\simeq
Z_{A_G}^\beta(\beta\uparrow\lambda),
\\
\Dist(G_1)\otimes_{\Dist(P_1^+)}T_{I,\lambda}^\omega
Z_{I,A_I}^\beta(\lambda)\otimes_{A_I}A_G^\beta
&\simeq
T_{\lambda}^\omega
Z_{A_G}^\beta(\lambda),
\\
\Dist(G_1)\otimes_{\Dist(P_1^+)}T_{I,\lambda}^\omega
Z_{I,A_I}^\beta(\beta\uparrow\lambda)\otimes_{A_I}A_G^\beta
&\simeq
T_{\lambda}^\omega
Z_{A_G}^\beta(\beta\uparrow\lambda),
\end{align*}
the commutative diagram
\[
\xymatrix{
\Ext_{\cC(\Omega_I,A_I^\beta)}^1
(Z_{I,A_I}^\beta(\lambda),
Z_{I,A_I}^\beta(\beta\uparrow\lambda)
\ar[rr]^-{t_{I,0}^\beta[\omega,\lambda]}
\ar[d]_-{\Dist(G_1)\otimes_{\Dist(P_1^+)}?\otimes_{A_I}A_G^\beta}&&
H_\beta^{-1}A_I^\beta/A_I^\beta
\ar[d]
\\
\Ext_{\cC(\Omega,A_G^\beta)}^1
(Z_{A_G}^\beta(\lambda),
Z_{A_G}^\beta(\beta\uparrow\lambda)
\ar[rr]^-{t_{0}^\beta[\omega,\lambda]}
&&
H_\beta^{-1}A_G^\beta/A_G^\beta
}
\]
follows.
More precisely,
if
$Y_{A_G}^\beta(\lambda)$
(resp.
$Y_{I,A_I}^\beta(\lambda)$)
is the module representing
$e^\beta(\lambda)$
(resp.
$e_I^\beta(\lambda)$),
then
$\Dist(G_1)\otimes_{\Dist(P_1^+)}Y_{I,A_I}^\beta(\lambda)
\otimes_{A_I}A_G^\beta=
Y_{A_G}^\beta(\lambda)$
with
$Y_{I,A_I}^\beta(\lambda)
\otimes_{A_I}A_G^\beta=
Y_{A_G}^\beta(\lambda)_\frn$
by (3.4).

\setcounter{equation}{0}
\noindent
(3.6)
We are now ready to show
\begin{thm}
There is a functorial isomorphism 
from $\cD(\Omega_I,A_I)$
to $\cD(\Omega,A_G)$
\[
\cV_{\Omega}\circ\hat\nabla_{P,A_I}
\simeq
\cI
\circ\cV_{\Omega_I}.
\]

\end{thm}

\pf
For each
$X\in\cD(\Omega, A_G)$
and $M\in\cD(\Omega_I, A_I)$
\begin{align*}
\cC(\Omega, A_G)
&
(X, \hat\nabla_{P,A_G}(M))
\simeq
\cC(\Omega, A_G)
(X, (\Dist(G_1)\otimes_{\Dist(P_1^+)}M^\tau)^\tau\otimes_{A_I}A_G)
\\
&
\simeq
\cC(\Omega, A_G)
(X, \{
\Dist(G_1)\otimes_{\Dist(P_1^+)}(M\otimes_{A_I}A_G)^\tau\}^\tau)
\\
&\simeq
\cC(\Omega, A_G)
(\Dist(G_1)\otimes_{\Dist(P_1^+)}(M\otimes_{A_I}A_G)^\tau, X^\tau)
\quad\text{by \cite[4.5.5]{AJS}}
\\
&\simeq
\cC_{L_I}(A_G)
((M\otimes_{A_I}A_G)^\tau, \Ann_{X^\tau}(\Dist^+(U_1^+(w_I)))),
\end{align*}
with
\begin{align*}
\Ann_{X^\tau}
(\Dist^+(U_1^+(w_I)))
&=
\{f\in
\Mod_{A_G}(X,A_G)
\mid
0=xf=f(\tau(x)?)
\ \forall x\in
\Dist^+(U_1^+(w_I))\}
\\
&\simeq
\{
X/\tau(\Dist^+(U_1^+(w_I)))X
\}^\tau
\\
&=
\{
X/(\Dist^+(U_1(w_I))X)
\}^\tau
=
(X_\frn)^\tau.
\end{align*}
Thus
\begin{align}
\cC(\Omega, A_G)
&
(X, \hat\nabla_{P,A_I}(M))
\simeq
\cC(\Omega_I, A_G)
((M\otimes_{A_I}A_G)^\tau,
(X_\frn)^\tau)
\\
\notag
&\simeq
\cC(\Omega_I, A_G)
(X_\frn,
M\otimes_{A_I}A_G).
\end{align}
It follows for each
$\lambda\in\Omega$
that
\begin{align*}
(\cV_{\Omega}\circ\hat\nabla_{P,A_I})
&
(M)(\lambda)
=
\cC(\Omega, A_G^\emptyset)(Z^\emptyset(\lambda),
\hat\nabla_{P,A_I}(M)^\emptyset)
\\
&\simeq
\cC(\Omega_I, A_G^\emptyset)(Z_{I,A_G^\emptyset}(\lambda),
M^\emptyset\otimes_{A_I}A_G^\emptyset)
\\
&\simeq
\cC(\Omega_I, A_I^{\emptyset})(Z_{I,A_I}^\emptyset(\lambda),
M^\emptyset)\otimes_{A_I}A_G^\emptyset
\quad\text{by \cite[3.2]{AJS}}
\\
&=
\begin{cases}
\cV_{\Omega_I}(M)(\lambda)\otimes_{A_I}A_G^\emptyset
&\text{if $\lambda\in\Omega_I$}
\\
0
&\text{else}
\end{cases}
\\
&=
(\cI
\circ\cV_{\Omega_I})(M)(\lambda).
\end{align*}

Assume $\lambda\in\Omega_I$.
If
$\beta\in R_I^+$,
\begin{align*}
(\cV_{\Omega}\circ
&
\hat\nabla_{P,A_I})
(M)(\lambda, \beta)
=
\cC(\Omega, A_G^\beta)(Y_{A_G}^\beta(\lambda),
\hat\nabla_{P,A_I}
(M)^\beta)
\\
&\simeq
\cC(\Omega_I, A_I^{\beta})(Y_{I,A_I}^\beta(\lambda),
M^\beta)\otimes_{A_I}A_G^\beta
\quad\text{likewise 
by 
 (3.5)
}
\\
&\simeq
\cV_{\Omega_I}
(M)(\lambda, \beta)\otimes_{A_I}A_G^\beta
=
(\cI
\circ\cV_{\Omega_I})(M)(\lambda, \beta).
\end{align*}
If $\beta\in R^+\setminus R_I^+$,
\begin{align*}
(\cV_{\Omega}\circ
&
\hat\nabla_{P,A_I})
(M)(\lambda, \beta)
=
\cC(\Omega, A_G^\beta)(Y_{A_G}^\beta(\lambda),
\hat\nabla_{P,A_I}
(M)^\beta)
\\
&\simeq
\cC(\Omega_I, A_G^\beta)((M^\beta\otimes_{A_I}A_G^\beta)^\tau,
(Y_{A_G}^\beta(\lambda)_\frn)^\tau)
\\
&\simeq
\cC(\Omega_I, A_G^\beta)((M^\emptyset\otimes_{A_I}A_G^\beta)^\tau,
(Z_{I,A_I}^\emptyset(\lambda)\otimes_{A_I}A_G^\beta)
^\tau)
\quad\text{as
$A_G^\beta\simeq
A_I^{\emptyset}\otimes_{A_I}A_G^\beta$
}
 \\
&\hspace{1cm}
\text{ in this case,
and hence
$Y_{A_G}^\beta(\lambda)_\frn\simeq
(Z_{I,A_I}^\emptyset(\lambda)
\oplus
Z_{I,A_I}^\emptyset(\beta\uparrow\lambda))\otimes_{A_I}A_G^\beta$}
\\
&\simeq
\cC(\Omega_I, A_I^{\beta})(Z_{I,A_I}^\emptyset(\lambda),
M^\emptyset)\otimes_{A_I}A_G^\beta
\\
&=
\cV_{\Omega_I}(M)(\lambda)\otimes_{A_I}A_G^\beta
=
(\cI
\circ\cV_{\Omega_I})(M)(\lambda, \beta).
\end{align*}

If $\lambda\in\Omega\setminus\Omega_I$
and if
$\beta\in R^+\setminus R_I^+$
with
$\beta\uparrow\lambda\in\Omega_I$, we have likewise
\[
(\cV_{\Omega}\circ
\hat\nabla_{P,A_I})
(M)(\lambda, \beta)
\simeq
\cV_{\Omega_I}(M)(\beta\uparrow\lambda)\otimes_{A_I}A_G^\beta
=
(\cI
\circ\cV_{\Omega_I})(M)(\lambda, \beta).
\]
If
$\lambda\in\Omega\setminus\Omega_I$
and if
$\beta\uparrow\lambda\notin\Omega_I$, 
\[
(\cV_{\Omega}\circ
\hat\nabla_{P,A_I})
(M)(\lambda, \beta)
=0=
(\cI\circ\cV_{\Omega_I})(M)(\lambda, \beta).
\]

\setcounter{equation}{0}
\noindent
(3.7)
Define
finally
a functor
$\tilde\cI
: \tilde\cK(\Omega_I,
S_{I,\Bbbk})\to
\tilde\cK(\Omega, S_\Bbbk)$
just like $\cI$
as follows:
for each
$\cM\in\cK(\Omega_I, A_I)$ 
and $\lambda\in\Omega$
set 
\[
(\tilde\cI
(\cM))(\lambda)
=
\begin{cases}
\cM(\lambda)\otimes_{S_{I,\Bbbk}
}
S_\Bbbk^\emptyset
&\text{if
$\lambda\in\Omega_I$}
\\
0
&\text{else},
\end{cases}
\]
and for each $\beta\in R^+$
set
\[
(\tilde\cI
(\cM))(\lambda,\beta)
=
\begin{cases}
\cM(\lambda,\beta)\otimes_{S_{I,\Bbbk}
}
S_\Bbbk^\beta
&\text{if
$\lambda\in\Omega_I$
and $\beta\in R_I^+$}
\\
\cM(\lambda)\otimes_{S_{I,\Bbbk}}
S_\Bbbk^\beta
&\text{if
$\lambda\in\Omega_I$
and $\beta\notin R_I^+$}
\\
\cM(\beta\uparrow\lambda)\otimes_{S_{I,\Bbbk
}}
S_\Bbbk^\beta
&\text{if
$\beta\uparrow\lambda\in\Omega_I$
and $\beta\notin R_I^+$}
\\
0
&\text{else}.
\end{cases}
\]
From
\cite[14.10]{AJS}
one has, in particular, for each
$\lambda\in\Omega_I$
\begin{equation}
\tilde\cI
(\cZ_{I,\lambda}^{w_I})
\simeq
\cZ_{\lambda}^{w_0}.
\end{equation}

Let
$Q_{I,A_I}=\bigoplus_{w\in W_I}Q_{I,A_I}(w\bullet\lambda_I^+)
\in\cC(\Omega_I, A_I)$
with
$Q_{I,A_I}(w\bullet\lambda^+)$
the lift 
of the projective cover of
$\hat L^P(w\bullet\lambda^+)$
for $L_I$
over $A_I$.
Let
$\cP_I$ be a graded
$S_{I,\Bbbk}$-form of
$\cV_{\Omega_I}(Q_{I,A_I})
$.

\begin{lem}
One has  
an isomorphism in $\cK(\Omega, A_G)$
\[
\tilde\cI
(\cP_I)\otimes_{S_\Bbbk}A_G
\simeq
\cI
(\cV_{\Omega_I}(Q_{I,A_I}))\simeq
\cV_{\Omega}\circ\hat\nabla_{P,A_I}(Q_{I,A_I}).
\]

\end{lem}

\pf
The first isomorphism follows from the definition
that
$\cP_I\otimes_{S_{I,\Bbbk}}A_I\simeq\cV_{\Omega_I}(Q_{I,A_I})$, and the
second from (3.6).

\setcounter{equation}{0}
\noindent
(3.8)
Let $E_{\Omega_I}=
\tilde\cK(\Omega_I, S_{I,\Bbbk})^\sharp(\cP_I, \cP_I)^\op$,
which is a $(Y_I\times\bbZ)$-graded
$S_{I,\Bbbk}$-algebra of finite type,
responsible for the structure of
$(Y_I\times\bbZ)$-graded $\Bbbk$-algebra
on
$E_{\Omega_I,\Bbbk}\simeq
E_{\Omega_I}\otimes_{S_{I,\Bbbk}}\Bbbk$.
Now set
$J=\tilde\cK(\Omega, S_{\Bbbk})^\sharp(\cP, \tilde\cI(\cP_I))$,
which comes equipped with a structure of
$(Y\times\bbZ)$-graded left $E_\Omega$
and
$(Y_I\times\bbZ)$-graded right $E_{\Omega,I}$-bimodule.
If
$J_\Bbbk=J\otimes_{S_\Bbbk}\Bbbk$,
it is 
thus
a $(Y\times\bbZ)$-graded left
$E_{\Omega,\Bbbk}$
and 
$(Y_I\times\bbZ)$-graded
right $E_{\Omega_I,\Bbbk}$-bimodule.

\begin{cor}
The parabolic induction functor
$\hat\nabla_P$ is
$\bbZ$-graded by the
bimodule
$J_\Bbbk$
in such a way
\[
\xymatrix{
\tilde\cC_\Bbbk(\Omega_I)
\ar[rrr]^-{J_\Bbbk\otimes_{E_{\Omega_I,\Bbbk}}?}
\ar[d]_-{\bar v_I}
\ar
@{}
[ddrrr]
|{\circlearrowright}
&&&
\tilde\cC_\Bbbk(\Omega)
\ar[dd]^-{\bar v}
\\
L_{I,1}T\Mod
\ar[d]
&&&
\\
P_{1}T\Mod
\ar[rrr]_{\hat\nabla_P}&&&
G_{1}T\Mod
}
\]
that
for each
$\lambda\in\Omega_I$
there is an isomorphism in
$\tilde\cC_\Bbbk(\Omega)$
\[
J_\Bbbk\otimes_{E_{\Omega_I,\Bbbk}}\tilde\nabla_{I,\Bbbk}(\lambda)
\simeq
\tilde\nabla_\Bbbk(\lambda)\langle\delta(\lambda)-\delta_I(\lambda)\rangle.
\]

\end{cor}

\pf
The commutativity of the diagram follows from
(3.2) by the isomorphism of
left
$E_{\Omega,\Bbbk}$
and right $E_{\Omega_I,\Bbbk}$-bimodules 
\begin{align*}
J_\Bbbk
&\simeq
J\otimes_{S_\Bbbk}A_G\otimes_{A_G}\Bbbk
\\
&\simeq
\cK(\Omega, A_G)^\sharp(\cV_\Omega(Q_{A_G}), \tilde\cI_\Omega(\cP_I)\otimes_{S_\Bbbk}A_G))
\otimes_{A_G}\Bbbk
\quad\text{by \cite[18.9.3]{AJS}}
\\
&\simeq
\cK(\Omega, A_G)^\sharp(\cV_\Omega(Q_{A_G}), \cV_\Omega\circ\hat\nabla_{P,A_I}(Q_{I,A_I}))
\otimes_{A_G}\Bbbk
\quad\text{by (3.7)}
\\
&\simeq
\cC(\Omega, A_G)^\sharp(Q_{A_G}, \hat\nabla_{P,A_I}(Q_{I,A_I}))
\otimes_{A_G}\Bbbk
\quad\text{by \cite[18.9.5/6]{AJS}}
\\
&\simeq
\cC_\Bbbk(\Omega)^\sharp(Q_{A_G}\otimes_{A_G}\Bbbk, \hat\nabla_{P,A_I}(Q_{I,A_I})\otimes_{A_G}\Bbbk)
\quad\text{by \cite[3.3]{AJS}}
\\
&\simeq
\cC_\Bbbk(\Omega)^\sharp(Q, \hat\nabla_{P}(Q_{I}))
\quad\text{
from (3.3.3)}.
\end{align*}
Also,
\begin{align*}
J_\Bbbk
\otimes_{E_{\Omega_I,\Bbbk}}
&
\tilde\nabla_{I,\Bbbk}(\lambda)
\\
&\simeq
(J
\otimes_{E_{\Omega_I}}
\tilde\nabla_{I,S_{I,\Bbbk}}(\lambda))\otimes_{S_{I,\Bbbk}}\Bbbk
=
\{
\tilde\cK(\Omega, S_\Bbbk)^\sharp(\cP, \tilde\cI_\Omega(\cP_I))\otimes_{E_{\Omega_I}}
\tilde\nabla_{I,S_{I,\Bbbk}}(\lambda)\}\otimes_{S_{I,\Bbbk}}\Bbbk
\\
&=
\{
\tilde\cK(\Omega, S_\Bbbk)^\sharp(\cP, \tilde\cI(\cP_I))\otimes_{E_{\Omega_I}}
\tilde\cK(\Omega_I, S_{I,\Bbbk})^\sharp(\cP_I, \tilde\cZ_{I,\lambda}^{w_I}\langle-\delta_I(\lambda)\rangle)\}
\otimes_{S_{I,\Bbbk}}\Bbbk
\\
&=
\tilde\cK(\Omega, S_\Bbbk)^\sharp(\cP, \tilde\cI_\Omega(\tilde\cZ_{I,\lambda}^{w_I}\langle-\delta_I(\lambda)\rangle)
\otimes_{S_{\Bbbk}}\Bbbk
\quad\text{by the five lemma}
\\
&=
\tilde\cK(\Omega, S_\Bbbk)^\sharp(\cP, \tilde\cZ_{\lambda}^{w_0}\langle-\delta_I(\lambda)\rangle)
\otimes_{S_{\Bbbk}}\Bbbk
\quad\text{by (3.7.1)}
\\
&=
\tilde\nabla_{S_\Bbbk}(\lambda)\langle\delta(\lambda)-\delta_I(\lambda)\rangle
\otimes_{S_{\Bbbk}}\Bbbk
=
\tilde\nabla_\Bbbk(\lambda)\langle\delta(\lambda)-\delta_I(\lambda)\rangle.
\end{align*}

\begin{center}
$4^\circ$
{\bf Rigidity}
\end{center}

Keep the notations of
\S3.
We will show that
all 
$\hat\nabla_P(\hat L^P(
\lambda))$
for $p$-regular
$\lambda\in\Lambda$
are $\bbZ$-graded.
Assuming Lusztig's conjecture on the irreducible character formulae for
$G_1T$-modules
\cite{L80}/\cite{Kato}, \cite{AJS}
has shown that
the endomorphism algebra
of a projective
$
Y$-generator 
for the block of $\lambda$ is
Koszul.
We show that the rigidity of $\hat\nabla_P(\hat L^P(
\lambda))$ 
follows from a result of
\cite{BGS}.  
The Lusztig conjecture is now a theorem
for large $p$
thanks to
\cite{AJS},
\cite{KL},
\cite{KT},
\cite{L94}
and more recently \cite{F}.
We will also 
find the Loewy length of
$\hat\nabla_P(\hat L^P(\lambda))$ for a $p$-regular $\lambda\in\Lambda$
to be uniformly
$\ell(w^I)+1$.

Thus fix a $p$-regular weight $\lambda$
and put
$\Omega=W_p\bullet\lambda$.
For
$M\in\cC_\Bbbk(\Omega)$
we let
$[M:\hat L(\mu)]$,
$\mu\in\Omega$,
denote the multiplicity of
simple $\hat L(\mu)$ among the 
$\cC_\Bbbk(\Omega)$-composition factors of $M$.

\setcounter{equation}{0}
\noindent
(4.1)
Let us first recall the construction of
$\tilde L_\Bbbk(\lambda)$,
which is a slight simplication over the one in \cite[18.12]{AJS}.
As $\tilde
Z_\Bbbk(\lambda)$ has a simple head in $E_{\Omega,\Bbbk}\bfmodgr_Y$
by the categorical equivalence
$v$,
the radical $\rad_{E_{\Omega,\Bbbk}\bfmodgr_Y}\tilde Z_\Bbbk(\lambda)$
of
$\tilde Z_\Bbbk(\lambda)$
in the category
$E_{\Omega,\Bbbk}\bfmodgr_Y$
is maximal.
But
$\rad_{E_{\Omega,\Bbbk}\bfmodgr_Y}\tilde Z_\Bbbk(\lambda)$
belongs to
$\tilde\cC_\Bbbk(\Omega)$
by \cite[E.11]{AJS},
and
hence coincides with the
radical $\rad_{\tilde\cC_\Bbbk(\Omega)}\tilde
Z_\Bbbk(\lambda)$
of
$\tilde
Z_\Bbbk(\lambda)$
in the category of
$\tilde\cC_\Bbbk(\Omega)$.
We set
$\tilde
L_\Bbbk(\lambda)=\tilde Z_\Bbbk(\lambda)/
\rad_{\tilde\cC_\Bbbk(\Omega)}\tilde
Z_\Bbbk(\lambda)$.
Then
$\bar v\tilde L_\Bbbk(\lambda)
=
Q\otimes_{E_{\Omega,\Bbbk}}\tilde L_\Bbbk(\lambda)
\simeq
\{
Q\otimes_{E_{\Omega,\Bbbk}}\tilde Z_\Bbbk(\lambda)\}
/\{
Q\otimes_{E_{\Omega,\Bbbk}}\rad_{E_{\Omega,\Bbbk}\bfmodgr_Y}
\tilde
Z_\Bbbk(\lambda)
\}
\simeq
Z_\Bbbk(\lambda)/\rad_{\cC_\Bbbk(\Omega)}
Z_\Bbbk(\lambda)
\simeq
\hat L(\lambda)
$.
In turn, 
$\tilde L_\Bbbk
(\lambda)
\simeq
H_{\Omega,\Bbbk}\hat L(\lambda)
$
in
$E_{\Omega,\Bbbk}\bfmodgr_Y$
while
$
H_{\Omega,\Bbbk}\hat L(\lambda)
=
\cC_\Bbbk(\Omega)^\sharp
(Q,\hat L(\lambda))
\simeq
\cC_\Bbbk(\Omega)(Q_\Bbbk(\lambda),\hat L(\lambda))$
as
$Q_\Bbbk(\lambda)$
is the projective cover of
$\hat L(\lambda)$,
and hence
$\tilde L_\Bbbk(\lambda)$
is of dimension 1.

By the equivalence
$v$ the $\tilde L_\Bbbk(\lambda)$,
$\lambda\in\Omega$,
exhaust the simple objects of
$E_{\Omega,\Bbbk}\bfmodgr_Y$.
If
$\tilde L$
is a simple object of
$\tilde\cC_\Bbbk(\Omega)$,
then 
\begin{align*}
0
&
\ne
E_{\Omega,\Bbbk}\bfmodgr_{Y}
(\tilde L,\tilde L_\Bbbk(\lambda))
\quad\text{for some
$\lambda\in\Omega$}
\\
&=\bigoplus_{i\in\bbZ}
E_{\Omega,\Bbbk}\bfmodgr_{Y}
(\tilde L,\tilde L_\Bbbk(\lambda))_i=
\bigoplus_{i\in\bbZ}
\tilde\cC_\Bbbk(\Omega)(\tilde L,\tilde L_\Bbbk(\lambda)\langle-i\rangle),
\end{align*}
and hence
$\tilde\cC_\Bbbk(\Omega)(\tilde L,\tilde L_\Bbbk(\lambda)\langle
i\rangle)\ne0$ for some
$i$,
Then
$\tilde L\simeq
\tilde L_\Bbbk(\lambda)\langle
i\rangle$
in $\tilde\cC_\Bbbk(\Omega)$
by their simplicity.
Such $\lambda$ and $i$ are unique
by
\cite[18.8]{AJS}.
Thus we have obtained the first 2 parts of

\begin{prop}
(i)
Each $\tilde L_\Bbbk(\lambda)$,
$\lambda\in\Omega$,
is 1-dimensional.

(ii)
Each simple object of
$\tilde\cC_\Bbbk(\Omega)$
is isomorphic to some 
$\tilde L(\lambda)\langle i\rangle$
for unique $\lambda\in\Omega$ and $i\in\bbZ$.
Any simple object of
$E_{\Omega,\Bbbk}\bfmodgr_{Y}$
is isomorphic to some 
$\tilde L(\lambda)$
for unique $\lambda\in\Omega$.

(iii) If
$M\in
\tilde\cC_\Bbbk(\Omega)$,
the radical
(resp. socle) series of
$M$ in $E_{\Omega,\Bbbk}\bfmodgr_{Y}$
and in $\tilde\cC_\Bbbk(\Omega)$ coincide.

\end{prop}

\pf
(iii)
We show first that
each radical layer
$\rad^i_{E_{\Omega,\Bbbk}\bfmodgr_{Y}}M/
\rad^{i+1}_{E_{\Omega,\Bbbk}\bfmodgr_{Y}}M$
remains semisimple in
$\tilde\cC_\Bbbk(\Omega)$.
As it
inherits the structure of
$\tilde\cC_\Bbbk(\Omega)$
from $M$
by
\cite[E.11]{AJS},
we may assume $M$ is semisimple in $E_{\Omega,\Bbbk}\bfmodgr_{Y}$.
If $L$ is a simple component of
$M$
in $E_{\Omega,\Bbbk}\bfmodgr_{Y}$,
as $L$ is 1-dimensional by (i),
each
$(E_{\Omega,\Bbbk})_{Y\times\{i\}}$,
$i\ne0$
annihilates $L$
while each element of $(E_{\Omega,\Bbbk})_{Y\times\{0\}}$ is acting by a scalar, and hence
$M$ is semisimple also in
$\tilde\cC_\Bbbk(\Omega)$;
each $\bbZ$-homogeneous component
$M_i$ of $M$ must be
$E_{\Omega,\Bbbk}$-stable.
On the other hand,
each
$\rad^i_{\tilde\cC_\Bbbk(\Omega)}M/
\rad^{i+1}_{\tilde\cC_\Bbbk(\Omega)}M$ is semisimple
in 
$E_{\Omega,\Bbbk}\bfmodgr_{Y}$
as each simple component is 1-dimensional by (i) again.
It now follows that the radical series of $M$ in
$E_{\Omega,\Bbbk}\bfmodgr_{Y}$ and
$\tilde\cC_\Bbbk(\Omega)$
coincide.

The socle version of
\cite[E.11]{AJS} holds, and hence
also the assertion about the socle series of $M$.

\setcounter{equation}{0}
\noindent
(4.2)
Assume now Lusztig's conjecture on the irreducible characters of
$G_1T$-modules.
Then $E_{\Omega,\Bbbk}$ is Koszul with respect to its $\bbZ$-gradation
thanks to 
\cite[18.17]{AJS}.
In particular,
$E_{\Omega,\Bbbk}$ is positively graded:
$E_{\Omega,\Bbbk}=\bigoplus_{i\in\bbN}(E_{\Omega,\Bbbk})_i$
with
$(E_{\Omega,\Bbbk})_0=\prod_{w\in W}\Bbbk\pi_w$,
and 
is generated by
$(E_{\Omega,\Bbbk})_1$ over $\Bbbk$ by \cite[Props. 2.1.3 and 2.3.1]{BGS}, where
$\pi_w:
\prod_{x\in W}
Q_\Bbbk(x\bullet
\lambda^+)\to
Q_\Bbbk(w\bullet
\lambda^+)$
is the projection.
Let $E_{\Omega,\Bbbk}\bfmodgr_\bbZ$ denote the category of finite dimensional
$\bbZ$-graded
$E_{\Omega,\Bbbk}$-modules.

\begin{prop}
Assume the Lusztig conjecture.

(i)
Each $\tilde L_\Bbbk(\lambda)$,
$\lambda\in\Omega$,
is homogeneous  of degree 0
with respect to the
$\bbZ$-grading.
In particular,
each $\tilde L_\Bbbk(w\bullet
\lambda^+)$,
$w\in W$,
is isomorphic to
$\Bbbk\pi_w$
in
$E_{\Omega,\Bbbk}\bfmodgr_\bbZ$.

(ii)
Each simple object of
$E_{\Omega,\Bbbk}\bfmodgr_\bbZ$ is isomorphic to some
$\tilde L_\Bbbk(w\bullet
\lambda^+)\langle i\rangle$
for unique $w\in W$ and $i\in\bbZ$.

(iii)
If
$M\in
\tilde\cC_\Bbbk(\Omega)$,
the radical
(resp. socle) series of
$M$ in $E_{\Omega,\Bbbk}\bfmodgr_{\bbZ}$
and in $\tilde\cC_\Bbbk(\Omega)$ coincide.

\end{prop}

\pf
(i)
Recall from (4.1) that the
$\bbZ$-grading on
$E_{\Omega,\Bbbk}$ arises from that of
$E_\Omega$.
Thus
$\Bbbk\pi_w=\Bbbk(\tilde\pi_w\otimes1)$
if
$\tilde\pi_w:\prod_{x\in
W}\cQ(x\bullet
\lambda^+)\to
\cQ(w\bullet
\lambda^+)$ is the projection.
But
\begin{align*}
\tilde\cK
&
(\Omega,S_\Bbbk)^\sharp(\bigoplus_{x\in W}\cQ(x\bullet
\lambda^+),
\cQ(w\bullet
\lambda^+))
\\
&=
\bigoplus_{x\in W}\bigoplus_{\gamma\in Y}
\bigoplus_{i\in\bbZ}
\cK
(\Omega,S_\Bbbk)(\cQ(x\bullet
\lambda^+)[\gamma],
\cQ(w\bullet
\lambda^+))_i
\quad\text{by definition
\cite[E.1, E.3]{AJS}}
\\
&=
\bigoplus_{x\in W}\bigoplus_{\gamma\in Y}
\cK
(\Omega,S_\Bbbk)(\cQ(x\bullet
\lambda^+)[\gamma],
\cQ(w\bullet
\lambda^+))
\quad\text{by
\cite[E.1]{AJS}}
\\
&\simeq
\bigoplus_{x\in W}\bigoplus_{\gamma\in Y}
\cK
(\Omega,S_\Bbbk)(\cQ(x\bullet
\lambda^++\gamma),
\cQ(w\bullet
\lambda^+))
\quad\text{by
\cite[17.6/18.5]{AJS}}
\end{align*}
with 
$\cK
(\Omega,S_\Bbbk)(\cQ(x\bullet
\lambda^++\gamma),
\cQ(w\bullet
\lambda^+))=
\bigoplus_{i>0}
\cK
(\Omega,S_\Bbbk)(\cQ(x\bullet
\lambda^++\gamma),
\cQ(w\bullet
\lambda^+))_i$
unless
$x\bullet
\lambda^++\gamma=w\bullet
\lambda^+$
while
$\cK
(\Omega,S_\Bbbk)(\cQ(w\bullet
\lambda^+),
\cQ(w\bullet
\lambda^+))_0
=
S_\Bbbk
\id_{\cQ(w\bullet
\lambda^+)}$
\cite[17.9]{AJS}.
On the other hand,
\begin{align*}
\tilde
Z_\Bbbk
&
(w\bullet
\lambda^++\gamma)
=
\tilde\cK(\Omega,S_\Bbbk)^\sharp(\bigoplus_{x\in W}\cQ(x\bullet
\lambda^+),
\cZ_{w\bullet
\lambda^++\gamma}\langle-\delta(w\bullet
\lambda^++\gamma)\rangle)\otimes_{S_\Bbbk}\Bbbk
\\
&\hspace{7cm}
\text{by definition
\cite[18.10.1 and 18.12]{AJS}}
\\
&\simeq
\bigoplus_{x\in W}
\bigoplus_{\nu\in Y}
\cK(\Omega,S_\Bbbk)(\cQ(x\bullet
\lambda^++\nu),
\cZ_{w\bullet
\lambda^++\gamma}\langle-\delta(w\bullet
\lambda^++\gamma)\rangle)
\otimes_{S_\Bbbk}\Bbbk
\quad\text{as above}.
\end{align*}
Each
$\cK(\Omega,S_\Bbbk)(\cQ(x\bullet
\lambda^++\nu),
\cZ_{w\bullet
\lambda^++\gamma}\langle-\delta(w\bullet
\lambda^++\gamma)\rangle)$
is a direct summand of
$\tilde\cK
(\Omega,S_\Bbbk)(\cQ(x\bullet
\lambda^++\nu),
\cQ(w\bullet
\lambda^++\gamma))$
by \cite[15.10 and 17.6.1/18.9.c]{AJS},
and hence
$\cK(\Omega,S_\Bbbk)(\cQ(x\bullet
\lambda^++\nu),
\cZ_{w\bullet
\lambda^++\gamma}\langle-\delta(w\bullet
\lambda^++\gamma)\rangle)=
\cK(\Omega,S_\Bbbk)(\cQ(x\bullet
\lambda^++\nu),
\cZ_{w\bullet
\lambda^++\gamma}\langle-\delta(w\bullet
\lambda^++\gamma)\rangle)_{>0}$
unless
$x\bullet
\lambda^++\nu=w\bullet
\lambda^++\gamma$,
i.e., $x=w$ and $\nu=\gamma$,
by \cite[17.9]{AJS} again
while
\begin{align*}
\cK(\Omega,S_\Bbbk)
&
(\cQ(w\bullet
\lambda^++\gamma),
\cZ_{w\bullet
\lambda^++\gamma}\langle-\delta(w\bullet
\lambda^++\gamma)\rangle)_0
\\
&\simeq
\cK(\Omega,S_\Bbbk)(\cZ'_{w\bullet
\lambda^++\gamma}\langle2|R^+|-\delta(w\bullet
\lambda^++\gamma)\rangle,
\cZ_{w\bullet
\lambda^++\gamma}\langle-\delta(w\bullet
\lambda^++\gamma)\rangle)_0
\\
&\hspace{8cm}
\text{by 
\cite[15.10 and 17.6.2]{AJS}}
\\
&\simeq
(S_\Bbbk)_0
\quad\text{by \cite[15.10.2]{AJS}}.
\end{align*}
Thus the epi
$\tilde
Z_\Bbbk(w\bullet
\lambda^++\gamma)/\tilde
Z_\Bbbk(w\bullet
\lambda^++\gamma)_{>0}\to
\tilde L_\Bbbk(w\bullet
\lambda^++\gamma)$
is an isomorphism of
$E_{\Omega,\Bbbk}\bfmodgr_{\bbZ}$
by dimension,
and hence $\tilde L_\Bbbk(w\bullet
\lambda^++\gamma)$ is of degree 0.
In particular,
$\tilde L_\Bbbk(w\bullet
\lambda^+)\simeq
\Bbbk(\tilde\pi_w\otimes1)$.

(ii)
Let $L$ be a simple object
of $E_{\Omega,\Bbbk}\bfmodgr_{\bbZ}$.
As $(E_{\Omega,\Bbbk})_{>0}L=0$, $L$ is a $(E_{\Omega,\Bbbk})_0$-module.
Then 
$L$ is
by its simplicity isomorphic to
some $\Bbbk\pi_w\langle i\rangle$, $w\in W$, $i\in\bbZ$.

(iii) now follows from (ii) just like 
(4.1.iii), applying \cite[E.11]{AJS}
to the pair
$(Y\times \bbZ, \bbZ)$
in place of
$(Y\times \bbZ, Y)$.

\setcounter{equation}{0}
\noindent
(4.3)
We are now to obtain 
from \cite[Prop. 2.4.1]{BGS}
the rigidity of 
$\hat\nabla_P(\hat L^P(
\lambda))$,
as well as
$\hat \nabla(\lambda)$
and $\hat Q(\lambda)=Q_\Bbbk(\lambda)$
for each
$\lambda\in\Omega$
demonstrated first in
\cite{AK89}
by a different method
using Vogan's version of the Lusztig conjecture.

\begin{lem}
Assume the Lusztig conjecture
on the irreducible characters of
$G_1T$-modules.
Let
$M\in
\tilde\cC_\Bbbk(\Omega)$.
If $M$ has a simple socle and a simple head
as an object of
$E_{\Omega,\Bbbk}\bfmodgr_{Y}$,
then
$M$ is rigid in
$E_{\Omega,\Bbbk}\bfmodgr_{Y}$.

\end{lem}

\pf
By the hypothesis
$M$ has a simple socle and a simple head
in
$(E_{\Omega,\Bbbk})\bfmodgr_\bbZ$ by (4.1) and (4.2).
If $\hd_{E_{\Omega,\Bbbk}\bfmodgr_\bbZ}M$
(resp. $\soc_{E_{\Omega,\Bbbk}\bfmodgr_\bbZ}M$) is concentrated in
degree $j$ (resp. $k$),
from
\cite[Prop.2.4.1]{BGS}
\[
\rad^i_{E_{\Omega,\Bbbk}\bfmodgr_\bbZ}M=M_{\geq i+j}
\quad\text{ and }\quad
\soc^i_{E_{\Omega,\Bbbk}\bfmodgr_\bbZ}M=M_{\geq k-i+1}
\quad\forall i.
\]
Thus
$M_{\geq k-i+1}
=\soc^i_{E_{\Omega,\Bbbk}\bfmodgr_\bbZ}M\geq
\rad^{\ell\ell(M)-i}_{E_{\Omega,\Bbbk}\bfmodgr_\bbZ}
M=M_{\geq \ell\ell(M)-i+j}$,
and hence
$k-i+1\leq\ell\ell(M)-i+j$.
As the equality holds for $i=0$,
$k+1=\ell\ell(M)+j$.
Then
$\forall i$,
$k-i+1\leq\ell\ell(M)-i+j=\ell\ell(M)-i+k+1-\ell\ell(M)=k-i+1$,
and hence
\[
\soc^i_{E_{\Omega,\Bbbk}\bfmodgr_Y}M=
\soc^i_{E_{\Omega,\Bbbk}\bfmodgr_\bbZ}M
=
M_{\geq\ell\ell(M)+j-i}
=
\rad^{\ell\ell(M)-i}_{E_{\Omega,\Bbbk}\bfmodgr_\bbZ}M
=
\rad^{\ell\ell(M)-i}_{E_{\Omega,\Bbbk}\bfmodgr_Y}M.
\]

\setcounter{equation}{0}
\noindent
(4.4)
Recalling from (1.4) that
each $\hat\nabla_P(\hat L^P(
\lambda))$
has a simple socle and a simple head yields

\begin{thm}
Assume the Lusztig conjecture.
Each  
$\hat\nabla_P(\hat L^P(
\lambda))$
for
$p$-regular
$\lambda$
is rigid.

\end{thm}

\setcounter{equation}{0}
\noindent
(4.5)
To determine eventually the Loewy series of
$\hat\nabla_P(\hat L^P(
\lambda))$,
we have to compute its Loewy length.
As
$\ell\ell(\hat\nabla_P(\hat L^P(
\lambda)))=
\ell\ell(^{w^I}\hat\nabla_P(\hat L^P(
\lambda)))
$,
we will compute
$\ell\ell(^{w^I}\hat\nabla_P(\hat L^P(
\lambda)))$.

\begin{lem}
$\hd_{G_1T}(^{w^I}\hat\nabla_P(\hat L^P(
\lambda)))=
\hat L(w^I
\bullet\lambda)\otimes-p(w^I\bullet0)$.

\end{lem}

\pf
We may assume
$\lambda^1=0$.
By (1.4)
\begin{align*}
\hd_{G_1T}
&
(^{w^I}\hat\nabla_P(\hat L^P(
\lambda)))
=
{^{w^I}\hd}_{G_1T}(\hat\nabla_P(\hat L^P(
\lambda)))
\\
&=
{^{w^I}}\{
\hat L(w^I
\bullet\lambda)
\otimes
p(
-2\rho_P+w_0((-w_I)\bullet
\lambda)^1-((-w_I)\bullet
\lambda)^1
))\}
\\
&=
{^{w^I}}\{
L((w^I\bullet
\lambda)^0)\otimes
p((w^I\bullet
\lambda)^0
-2\rho_P+w_0((-w_I)\bullet
\lambda)^1-((-w_I)\bullet
\lambda)^1))\}
\\
&=
L((w^I\bullet
\lambda)^0)\otimes
pw^I\{(w^I\bullet
\lambda)^0
-2\rho_P+w_0((-w_I)\bullet
\lambda)^1-((-w_I)\bullet
\lambda)^1
\}
\end{align*}
while
$\hat L(w^I
\bullet\lambda)\otimes-p(w^I\bullet0)=
L((w^I\bullet\lambda)^0)\otimes
p\{(w^I\bullet
\lambda)^1-(w^I\bullet0)\}$.
Thus we are to show
\begin{equation}
(w^I\bullet
\lambda)^1-(w^I\bullet0)
=
w^I\{(w^I\bullet
\lambda)^1
-2\rho_P+w_0((-w_I)\bullet
\lambda)^1-((-w_I)\bullet
\lambda)^1
\}.
\end{equation}
Write
$w_I\bullet
\lambda=\mu^0+p\mu^1$
with
$\mu^0\in\Lambda_p$
and
$\mu^1\in\Lambda$.
Thus
$\mu^0$ is $p$-regular.
As
$w^I\bullet
\lambda=
w_0\bullet(\mu^0+p\mu^1)=w_0\bullet\mu^0+pw_0\mu^1$,
$(w^I\bullet
\lambda)^1=w_0\mu^1-\rho$.
Likewise, 
as
$(-w_I)\bullet
\lambda=(-1)\bullet(\mu^0+p\mu^1)=(-1)\bullet\mu^0-p\mu^1$,
$((-w_I)\bullet
\lambda)^1=-\mu^1-\rho$.
It follows that
the LHS of (1) is equal by (1.1)
to
\[
w_0\mu^1-\rho-w_0w_I\bullet0
=
w_0\mu^1-\rho-(w_0w_I\rho-\rho)
=
w_0\mu^1-\rho-w_02\rho_P
=
w_0(\mu^1+\rho-2\rho_P)
\]
while the RHS of (1) is equal to
\begin{align*}
w^I
&
\{
w_0\mu^1-\rho-2\rho_P+w_0(-\mu^1-\rho)-(-\mu^1-\rho)\}
=
w^I(-2\rho_P+\mu^1+\rho)
\\
&=
w_0(w_I\mu^1+w_I\rho-2\rho_P)
\quad\text{as
$w_I2\rho_P=2\rho_P$ 
}.
\end{align*}
Thus we are left to verify that
$\mu^1+\rho=w_I(\mu^1+\rho)$,
for which we have only to check  
$\langle\mu^1+\rho,\alpha^\vee\rangle=0$
$\forall\alpha\in I$.
But
\begin{align*}
]0,p[
+p\langle\mu^1+\rho,\alpha\rangle
&\ni
\langle\mu^0+\rho,\alpha\rangle+p\langle\mu^1+\rho,\alpha\rangle
\quad\text{as
$\mu^0$
is $p$-regular}
\\
&=
\langle\mu^0+p\mu^1+p\rho+\rho,\alpha^\vee\rangle=
\langle
w_I\bullet
\lambda+p\rho+\rho,\alpha^\vee\rangle
=
\langle
\lambda+\rho,
w_I\alpha^\vee\rangle
+p
\\
&\in
-]0,p[+p
\quad\text{as
$w_I\alpha\in-I$}
\\
&=
]0,p[,
\end{align*}
and hence
$\langle\mu^1+\rho,\alpha^\vee\rangle=0$,
as desired.

\setcounter{equation}{0}
\noindent
(4.6)
Recall from (1.6.3) that
$^{w^I}\hat\nabla_P(\hat L^P(
\lambda))\leq
\hat\nabla_{w^I}((w^I\bullet
\lambda)
\langle
w^I\rangle)\otimes-p(w^I\bullet0)$.
Recall also from
\cite{AK89}
an intertwining homomorphism
$\phi_{w}\in
G_1T\Mod(\hat\nabla_{w}((w\bullet
\lambda)\langle w\rangle),
\hat\nabla(w\bullet
\lambda))\setminus0$
for each
$w\in W$,
which is unique up to
$\Bbbk^\times$.
As
$1=[\hat\nabla(w^I\bullet
\lambda):\hat
L(w^I\bullet
\lambda)]
=
[\hat\nabla_{w^I}((w^I\bullet
\lambda)\langle
w^I\rangle):\hat
L(w^I\bullet
\lambda)]$
by \cite[1.2.3]{AK89},
one obtains from (4.5) a commutative diagram of
$G_1T$-modules
\begin{equation}
\xymatrix{
\hat\nabla_{w^I}((w^I\bullet
\lambda)\langle
w^I\rangle)\otimes-p(w^I\bullet0)
\ar[rr]^-{\phi_{w^I}\otimes-p(w^I\bullet0)}
&&
\hat\nabla(w^I\bullet
\lambda)\otimes-p(w^I\bullet0)
\\
^{w^I}\hat\nabla_P(\hat L^P(
\lambda))
\ar@{->>}[rr]
\ar@{^(->}!<0ex,2.5ex>;[u]
&&
\hat
L(w^I\bullet
\lambda)\otimes-p(w^I\bullet0).
\ar@{^(->}[u]
}
\end{equation}
As
$\phi_{w^I}(\soc^{\ell(w^I)}
\hat\nabla_{w^I}((w^I\bullet
\lambda)\langle
w^I\rangle))=0$
\cite{AK89},
we must have
\begin{equation}
\ell\ell(^{w^I}\hat\nabla_P(\hat L^P(
\lambda)))\geq
\ell(w^I)+1.
\end{equation}
On the other hand,
there is another intertwining homomorphism
$\phi'_{w^I}\in
G_1T\Mod(\hat\nabla_{w_0}((w^I\bullet
\lambda)\langle
w_0\rangle),
\hat\nabla_{w^I}((w^I\bullet
\lambda)\langle
w^I\rangle))\setminus0$.
As
\begin{align*}
\hd_{G_1T}
&
\hat\nabla_{w_0}((w^I\bullet
\lambda)\langle
w_0\rangle)
\otimes
-p(w^I\bullet0)
\\
&=
\hd_{G_1T}
\hat\Delta(w^I\bullet
\lambda)
\otimes
-p(w^I\bullet0)
\quad\text{by \cite[1.2]{AK89}
}
\\
&=
\hat L(w^I\bullet
\lambda)
\otimes
-p(w^I\bullet0)
=
\hd_{G_1T}(^{w^I}\hat\nabla_P(\hat L^P(
\lambda))),
\end{align*}
one obtains as in (1) another commutative diagram
\begin{equation}
\xymatrix
{
\hat\nabla_{w_0}((w^I\bullet
\lambda)\langle
w_0\rangle)\otimes-p(w^I\bullet0)
\ar[dd]_-{\phi'_{w^I}\otimes-p(w^I\bullet0)}
\ar@{->>}[rd]
\\
&
^{w^I}\hat\nabla_P(\hat L^P(
\lambda))
\ar@{_(->}!<2ex,0ex>;[ld]
\\
\hat\nabla_{w^I}((w^I\bullet
\lambda)\langle
w^I\rangle)\otimes-p(w^I\bullet0).
}
\end{equation}
with
$\phi'_{w^I}(\soc^{\ell(w_0)-\ell(w^I)}
\hat\nabla_{w_0}((w^I\bullet
\lambda)\langle
w_0\rangle))=0$.
Assuming the Lusztig
conjecture we have 
$\ell\ell(\hat\nabla_{w_0}((w^I\bullet
\lambda)\langle
w_0\rangle))=\ell(w_0)+1$.
It follows
that
\[
\ell\ell(^{w^I}\hat\nabla_P(\hat L^P(
\lambda)))
\leq\ell(w_0)+1-\{\ell(w_0)-\ell(w^I)\}
=
\ell(w^I)+1.
\]
Thus, together with (2), we have obtained

\begin{thm}
Assume the Lusztig
conjecture. For any $p$-regular $\lambda\in\Lambda$
\[
\ell\ell(\hat\nabla_P(\hat L^P(\lambda)))=
\ell(w^I)+1.
\]

\end{thm}

\setcounter{equation}{0}
\noindent
(4.7)
{\bf Remark:}
This is a generalization of
\cite[1.4]{KY} and \cite{K09},
where we found for
$G$ of rank at most 2
or in case
$G=GL_{n+1}(\Bbbk)$ 
with $P$ a maximal parabolic such that
$G/P\simeq\bbP^n$ for any
$n\in\bbN$ that
$\ell\ell(\hat\nabla_P(\lambda))=\ell(w^I)+1$
for
$p$-regular
$\lambda\in\Lambda_P$. 
In fact, for $G/P\simeq\bbP^n$
we computed $\ell\ell(\hat\nabla_P(\lambda))$ for any
$\lambda\in\Lambda_P$
in
\cite[2.3]{K09}
dispensing with the Lusztig
conjecture.

\setcounter{equation}{0}
\noindent
(4.8)
Recall that
$\tilde Z_\Bbbk(\lambda)/\tilde Z_\Bbbk(\lambda)_{>0}\simeq
\tilde L_\Bbbk(\lambda)\simeq
\hd_{E_{\Omega,\Bbbk}\bfmodgr_\bbZ}\tilde Z_\Bbbk(\lambda)$
for each
$\lambda\in\Omega$.
It follows that
the $\bbZ$-gradation on
$\tilde Z_\Bbbk(\lambda)$ is such that
for each
$j\in\bbN
$
\begin{align*}
\tilde Z_\Bbbk(\lambda)_{\geq j}
&=
\rad_{E_{\Omega,\Bbbk}\bfmodgr_\bbZ}^j\tilde Z_\Bbbk(\lambda)
=
\rad_{E_{\Omega,\Bbbk}\bfmodgr_Y}^j\tilde
Z_\Bbbk(\lambda)
\\
&=
\soc_{E_{\Omega,\Bbbk}\bfmodgr_\bbZ}^{|R^+|+1-j}
\tilde Z_\Bbbk(\lambda)
=
\soc_{E_{\Omega,\Bbbk}\bfmodgr_Y}^{|R^+|+1-j}
\tilde
Z_\Bbbk(\lambda),
\end{align*}
and hence
\[
\soc_{\cC_\Bbbk(\Omega)}^{|R^+|+1-j}
Z_\Bbbk(\lambda)
=
\bar v(\tilde Z_\Bbbk(\lambda)_{\geq j})
=
\rad_{\cC_\Bbbk(\Omega)}^{j}
Z_\Bbbk(\lambda).
\]
More generally,

\begin{prop}
Assume the Lusztig conjecture.
The $\bbZ$-gradation on
each
$\tilde Z_\Bbbk^w(\lambda)$,
$\lambda\in\Omega$, $w\in W$,
is such that for each $i\in\bbN$
\begin{align*}
\rad_{\cC_\Bbbk(\Omega)}^i
Z_\Bbbk^w(\lambda\langle w\rangle)
&=
\rad_{G_1T}^i
\hat\nabla_{ww_0}(\lambda\langle ww_0\rangle))
=
\bar v(\tilde Z_\Bbbk^w(\lambda)_{\geq-\ell(w)+i})
\\
&=
\soc_{\cC_\Bbbk(\Omega)}^{|R^+|+1-i}
Z_\Bbbk^w(\lambda\langle w\rangle)
=
\soc_{G_1T}^{|R^+|+1-i}
\hat\nabla_{ww_0}(\lambda\langle ww_0\rangle).
\end{align*}
Thus $\forall\mu\in\Omega$,
\begin{align*}
[\rad_{\cC_\Bbbk(\Omega),i}Z_\Bbbk^w(\lambda\langle w\rangle):\hat L(\mu)]
&=
[\tilde Z_\Bbbk^w(\lambda):\tilde L_\Bbbk(\mu)\langle-\ell(w)+i\rangle]
\\
&=
[\soc_{\cC_\Bbbk(\Omega),|R^+|+1-i}Z_\Bbbk^w(\lambda\langle w\rangle):\hat L(\mu)],
\end{align*}
where the middle term is the multiplicity of
simple $\tilde L_\Bbbk(\mu)\langle-\ell(w)+i\rangle$
in
$\tilde Z_\Bbbk^w(\lambda)$
considered as objects of
$E_{\Omega,\Bbbk}\bfmodgr_\bbZ$.

\end{prop}

\pf
One has from
\cite[15.3.2]{AJS}
\[
\Bbbk=
\cC_\Bbbk(\Omega)
(Z_\Bbbk^w(\lambda\langle w\rangle),
Z_\Bbbk(\lambda))
=
\tilde\cC_\Bbbk(\Omega)
(\tilde
Z_\Bbbk^w(\lambda),
\tilde
Z_\Bbbk(\lambda)\langle-2\ell(w)\rangle).
\]
Let
$j\in\bbZ$ minimal such that
$\tilde
Z_\Bbbk^w(\lambda)_j\ne0$,
so
$\tilde
Z_\Bbbk^w(\lambda)_{\geq j}/
\tilde
Z_\Bbbk^w(\lambda)_{>j}=
\hd_{E_{\Omega,\Bbbk}\bfmodgr_\bbZ}\tilde
Z_\Bbbk^w(\lambda)
=
\hd_{E_{\Omega,\Bbbk}\bfmodgr_Y}\tilde
Z_\Bbbk^w(\lambda)
=
H_{\Omega,\Bbbk}(\rad_{\cC_\Bbbk(\Omega),0}Z_\Bbbk^w(\lambda\langle w\rangle))$,
which is sent to
\begin{align*}
\tilde
Z_\Bbbk(\lambda)_{\geq j+2\ell(w)}
&/
\tilde
Z_\Bbbk(\lambda)_{>j+2\ell(w)}
=
(\tilde
Z_\Bbbk(\lambda)\langle-2\ell(w)\rangle)_{\geq j}
/
(\tilde
Z_\Bbbk(\lambda)\langle-2\ell(w)\rangle)_{>j}
\\
&=
H_{\Omega,\Bbbk}
(\rad_{\cC_\Bbbk(\Omega),\ell(w)}Z_\Bbbk(\lambda))
\\
&=
\tilde Z_\Bbbk(\lambda)_{\geq\ell(w)}/
\tilde Z_\Bbbk(\lambda)_{>\ell(w)}
\quad\text{by above}.
\end{align*}
Thus $j=-\ell(w)$.
As
$\ell\ell(Z_\Bbbk^w(\lambda\langle w\rangle))=|R^+|+1$,
the assertion follows.

\setcounter{equation}{0}
\noindent
(4.9)
Untwisting
$w^I$ of (4.6.3) reads 
\begin{equation}
\xymatrix
{
Z_\Bbbk^{w_Iw_0}(
\lambda\langle
w_Iw_0\rangle
)
\ar[dd]
\ar@{->>}[rd]
\\
&
\hat\nabla_P(\hat L^P(
\lambda))
\ar@{_(->}!<2ex,1.5ex>;[ld]
\\
Z_\Bbbk^{w_0}(
\lambda\langle w_0\rangle
)).
}
\end{equation}
Thus one obtains a commutative diagram
in
$E_{\Omega,\Bbbk}\bfmodgr_Y$
\begin{equation}
\xymatrix
{
H_{\Omega,\Bbbk}
Z_\Bbbk^{w_Iw_0}(
\lambda\langle
w_Iw_0\rangle
)
\ar[dd]
\ar@{->>}[rd]
\\
&
H_{\Omega,\Bbbk}\hat\nabla_P(\hat L^P(
\lambda))
\ar@{_(->}!<0.8ex,0.5ex>;[ld]
\\
H_{\Omega,\Bbbk}
Z_\Bbbk^{w_0}(
\lambda\langle w_0\rangle
).
}
\end{equation}
Recall 
that
$\cC_\Bbbk(\Omega)(
Z_\Bbbk^{w_Iw_0}
(
\lambda\langle
w_Iw_0\rangle
),
Z_\Bbbk^{w_0}
(
\lambda\langle
w_0\rangle
))$
is 1-dimensional.
On the other hand, 
each
$Z_\Bbbk^w(\lambda\langle
w\rangle)$,
$w\in W$,
admits a graded object 
$\tilde Z_\Bbbk^w(\lambda)\in\tilde\cC_\Bbbk(\Omega)$ such that
$\bar v\tilde Z^w_\Bbbk(\lambda)
\simeq
Z_\Bbbk^w(\lambda\langle
w\rangle)$.
It follows that
\begin{align*}
E_{\Omega,\Bbbk}\bfmodgr_Y
&(
H_{\Omega,\Bbbk}
Z_\Bbbk^{w_Iw_0}(
\lambda\langle
w_Iw_0\rangle
),
H_{\Omega,\Bbbk}
Z_\Bbbk^{w_0}(
\lambda\langle w_0\rangle
))
\\
&=
E_{\Omega,\Bbbk}\bfmodgr_Y
(
\tilde
Z^{w_Iw_0}_\Bbbk(
\lambda
),
\tilde
Z^{w_0}(
\lambda
))
\\
&=
\bigoplus_{i\in\bbZ}
E_{\Omega,\Bbbk}\bfmodgr_Y
(
\tilde
Z^{w_Iw_0}_\Bbbk(
\lambda
),
\tilde
Z^{w_0}(
\lambda
))_i
=
\bigoplus_{i\in\bbZ}
\tilde\cC_\Bbbk(\Omega)
(
\tilde
Z^{w_Iw_0}_\Bbbk(
\lambda
),
\tilde
Z^{w_0}(
\lambda
)\langle-i\rangle)
\\
&=
\tilde\cC_\Bbbk(\Omega)(
\tilde
Z^{w_Iw_0}_\Bbbk(
\lambda
),
\tilde
Z^{w_0}(
\lambda
)\langle
-j\rangle)
\end{align*}
for some single $j\in\bbZ$ by dimension;
in fact,
$j=0$
by \cite[15.3.2]{AJS}.
Then,
taking
$\eta\in
\tilde\cC_\Bbbk(\Omega)(
\tilde
Z^{w_Iw_0}_\Bbbk(
\lambda
),
\tilde
Z^{w_0}(
\lambda
))
\setminus0$,
$\im(\eta)\in\tilde\cC_\Bbbk(\Omega)$
with
$\bar v(\im(\eta))=\hat\nabla_P(\hat L^P(
\lambda))$.
This gives another proof that
$\hat\nabla_P(\hat L^P(
\lambda))$ is
$\bbZ$-graded, and hence is rigid.

\begin{cor}
Assume the Lusztig conjecture.
The $\bbZ$-gradation on
$\im(\eta)$
is such that for each $i\in\bbN$
\begin{align*}
\bar v((\im(\eta))_{\geq-i})
&=
\rad_{G_1T}^{\ell(w^I)-i}
\hat\nabla_P(\hat L^P(
\lambda))
=
\soc_{G_1T}^{i+1}
\hat\nabla_P(\hat L^P(
\lambda)).
\end{align*}

\end{cor}

\pf
As
\begin{align*}
\soc_{G_1T}
\hat\nabla_P(\hat L^P(
\lambda))
&=
\soc_{G_1T}
Z^{w_0}(
\lambda\langle w_0\rangle
)
\quad\text{by (1)}
\\
&=
\bar v(\tilde
Z^{w_0}(
\lambda\langle w_0\rangle
)_0)
\quad\text{by (4.8)},
\end{align*}
$\bar v(\im(\eta)_0)=
\soc_{G_1T}
\hat\nabla_P(\hat L^P(
\lambda))$,
and hence the assertion.

\begin{center}
$5^\circ$
{\bf The Loewy series}
\end{center} 

Keep the notations of \S4.
In this section we derive a formula to
describe the socle series of
$\hat\nabla_P(\hat L^P(
\lambda))$.
We continue to assume the Lusztig 
conjecture.

\setcounter{equation}{0}
\noindent
(5.1)
Let us first recall 
from \cite{AK89}
or from
\cite[18.19]{AJS} 
a formula for the socle series of
$\hat\nabla(
\lambda)$:
\begin{equation}
Q^{\mu,\lambda}
=
\sum_jq^{\frac{\rd(\mu, \lambda)-j}{2}}
[\soc_{j+1}\hat\nabla(\lambda) : \hat L(
\mu)],
\end{equation} 
where
$\rd(\mu, \lambda)=\rd(A,C)$
is the distance from alcove
$A$ containing $\mu$
to alcove
$C$ containing $\lambda$
and $Q^{\mu, \lambda}=Q^{A,C}$ is a periodic inverse Kazhdan-Lusztig polynomial defined 
in
\cite{L80}.

On the other hand,
the Lusztig conjecture for the $L_{I,1}T$-modules asserts
\begin{equation}
\ch\hat L^P(
\lambda)
=
\sum_{\mu\in W_{I,p}\bullet\lambda}
(-1)^{\rd_I(\mu,\lambda)}\hat P_{\mu,\lambda}^I(1)\ch\hat\nabla^P(\mu),
\end{equation}
where
$\rd_I(\mu,\lambda)$
is the distance from alcove
$A$ containing $\mu$
to alcove
$C$ containing $\lambda$
with respect to
$W_{I,p}$
and $\hat P^I_{\mu,\lambda}=\hat P^I_{A,C}$ is Kato's periodic Kazhdan-Lusztig polynomial 
for $W_{I,p}$
\cite{Kato}.
We will prove a formula
\begin{equation}
\sum_jq^{\frac{\rd(\mu, \lambda)-j}{2}}
[\soc_{j+1}\hat\nabla_P(\hat L^P(\lambda)) : \hat L(
\mu)]
=
\sum_{\nu\in
W_{I,p}\bullet\lambda}Q^{\mu,\nu}(-1)^{\rd_I(\nu,\lambda)}\hat P^I_{\nu,\lambda}.
\end{equation}
The formula 
reduces to (1) in case
$I=\emptyset$, i.e., 
when
$P$ is a Borel subgroup.
It also holds for $P=G$ by the inversion formula
$\sum_\nu
Q^{\mu,\nu}(-1)^{\rd(\nu,\lambda)}\hat P_{\nu,\lambda}=\delta_{\mu,\lambda}$
\cite[11.10]{L80}/\cite[p. 129]{Kato}.

\setcounter{equation}{0}
\noindent
(5.2)
Let us write
$\tilde\nabla_\Bbbk(\lambda)=\tilde
Z_\Bbbk^{w_0}(\lambda)$
for each 
$\lambda\in\Omega$.
Then
the formula (5.1.1) reads by (4.8)
with $q^{\frac{1}{2}}=t$
\begin{align}
Q^{\mu,\lambda}(t^2)
&=
\sum_j
t^{\rd(\mu,\lambda)-j}[\tilde\nabla_\Bbbk
(\lambda)
:\tilde L_\Bbbk(\mu)\langle
-j\rangle].
\end{align}
If we write
$Q^{\lambda,\nu}(t)=\sum_jQ^{\lambda\nu}_jt^{\frac{j}{2}}$ with
$Q^{\lambda\nu}_j\in\bbZ$,
the formula (1) reads
in the Grothendieck group of
$E_{\Omega,\Bbbk}\bfmodgr_\bbZ$
\begin{equation}
[\tilde\nabla_\Bbbk(\lambda)]=\sum_{j\in\bbZ}\sum_{\mu\in\Omega}Q^{\mu\lambda}_{\rd(\mu,\lambda)-j}[\tilde
L_\Bbbk(\mu)\langle-j\rangle],\end{equation}
inverting which reads,
if we write
$\hat P_{\mu,\lambda}(t)=
\sum_{j\in\bbZ}\hat P_{\mu\lambda, j}t^{\frac{1}{2}}$,
\begin{equation}
[\tilde L_\Bbbk(\lambda)]=\sum_{j\in\bbZ}\sum_{\mu\in\Omega}
(-1)^{\rd(\mu,\lambda)}\hat
P_{\mu\lambda,j+\rd(\mu, \lambda)}[\tilde\nabla_\Bbbk(\mu)\langle
 j\rangle].
\end{equation}

We now verify the formula (5.1.3).

\begin{thm}
Assume the Lusztig conjecture.
If
$\lambda$
is a $p$-regular weight, the Loewy series of
$\hat\nabla_P(\hat L^P(\lambda))$ is given by 
\[
\sum_{j\in\bbN}q^{\frac{\rd(\mu, \lambda)-j}{2}}
[\soc_{j+1}\hat\nabla_P(\hat L^P(\lambda)) : \hat L(
\mu)]
=
\sum_{\nu\in
W_{I,p}\bullet\lambda}Q^{\mu,\nu}(-1)^{\rd_I(\nu,\lambda)}\hat P^I_{\nu,\lambda}
\quad
\forall\mu\in\Lambda.
\]

\end{thm}

\pf
Put
$\tilde\nabla_P=
J_\Bbbk\otimes_{E_{\Omega_I,\Bbbk}}?$
from (3.8)
for simplicity.
In the Grothendieck group of
$\tilde\cC_\Bbbk(\Omega_I)$
the formula
(3) reads
$[\tilde
L_{I,\Bbbk}(\lambda)]=\sum_{j\in\bbZ,\mu\in\Omega}
(-1)^{\rd_I(\mu,\lambda)}\hat
P^I_{\mu\lambda,
 j+\rd_I(\mu,\lambda)}[\tilde
\nabla_{I,\Bbbk}(\mu)\langle
j\rangle]$.
Put $n_\lambda=\delta(\lambda)-
\delta_I(\lambda)$,
so
$\tilde\nabla_P(\tilde\nabla_{I,\Bbbk}(\lambda))\simeq
\tilde\nabla_\Bbbk(\lambda)\langle
n_\lambda\rangle$
by (3.8).
As
$\hat\nabla_P$ is exact,
so is
$\tilde\nabla_P$ by (3.8) also.
Then
\begin{align*}
[\tilde\nabla_P(\tilde
L_{I,\Bbbk}(\lambda))]
&=
\sum_{\mu,j}
(-1)^{\rd_I(\mu,\lambda)}\hat
P^I_{\mu\lambda,
 j+\rd_I(\mu,\lambda)}
 [\tilde\nabla_P(\tilde\nabla_{I,\Bbbk}(\mu))\langle j\rangle]
\\
&=
\sum_{\mu,j}
(-1)^{\rd_I(\mu,\lambda)}\hat
P^I_{\mu\lambda,
 j+\rd_I(\mu,\lambda)}
 [\tilde\nabla_\Bbbk(\mu)\langle n_\mu+j\rangle]
\\
&=
\sum_{\mu,j}
(-1)^{\rd_I(\mu,\lambda)}\hat
P^I_{\mu\lambda,
 j+\rd_I(\mu,\lambda)}
 \sum_{k,\nu}
Q^{\nu\mu}_{\rd(\nu,\mu)-k}[\tilde
L_\Bbbk(\nu)\langle n_\mu+j-k\rangle]
\\
&=
\sum_{\mu,j,\nu,k}
(-1)^{\rd_I(\mu,\lambda)}\hat
P^I_{\mu\lambda,
 j+\rd_I(\mu,\lambda)}
Q^{\nu\mu}_{\rd(\nu,\mu)-k}[\tilde
L_\Bbbk(\nu)\langle n_\mu+j-k\rangle].
\end{align*}
Recall now
$\im(\eta)$ from (4.9).
As
$\tilde L_{I,\Bbbk}(\lambda)\leq
\tilde\nabla_{I,\Bbbk}(\lambda)$,
$\tilde\nabla_P(\tilde L_{I,\Bbbk}(\lambda))\leq
\tilde\nabla_P(\tilde\nabla_{I,\Bbbk}(\lambda))
\simeq
\tilde\nabla_\Bbbk(\lambda)\langle n_\lambda\rangle$.
As $\im(\eta)\leq\tilde\nabla_\Bbbk(\lambda)$, it follows that
$\tilde\nabla_P(\tilde L_{I,\Bbbk}(\lambda))\simeq
\im(\eta)\langle n_\lambda\rangle$.
Thus
\begin{align*}
[\soc_{
i+1}\hat\nabla_P
&
(\hat L^P(\lambda)):\hat L(\nu)]
=
[\im(\eta):\tilde L_\Bbbk(\nu)\langle-i\rangle]
\quad\text{by (4.9)}
\\
&=
[\tilde\nabla_P(\tilde
L_{I,\Bbbk}(\lambda))\langle-n_\lambda\rangle:\tilde L_\Bbbk(\nu)\langle-i\rangle]
=
[\tilde\nabla_P(\tilde
L_{I,\Bbbk}(\lambda)):\tilde L_\Bbbk(\nu)\langle
n_\lambda-i\rangle]
\\
&=
\sum_{\mu,j}
(-1)^{\rd_I(\mu,\lambda)}\hat P^I_{\mu\lambda,j+\rd_I(\mu,\lambda)}
Q^{\nu\mu}_{\rd(\nu,\mu)+\delta(\lambda)-\delta_I(\lambda)-\delta(\mu)+\delta_I(\mu)-i-j}
\\
&=
\sum_{\mu,j}
(-1)^{\rd_I(\mu,\lambda)}\hat P^I_{\mu\lambda,j+\rd_I(\mu,\lambda)}
Q^{\nu\mu}_{\rd(\nu,\lambda)-\rd_I(\mu,\lambda)-i-j}
\\
&=
\sum_{\mu,j}
(-1)^{\rd_I(\mu,\lambda)}\hat P^I_{\mu\lambda
j}
Q^{\nu\mu}_{\rd(\nu,\lambda)-i-j},
\end{align*}
and hence
\begin{align*}
\sum_i
t^{\rd(\nu,\lambda)-i}
&
[\soc_{
i+1}\hat\nabla_P
(\hat L^P(\lambda)):\hat L(\nu)]
=
\sum_i
t^{\rd(\nu,\lambda)-i}
\sum_{\mu,j}
(-1)^{\rd_I(\mu,\lambda)}\hat P^I_{\mu\lambda
j}
Q^{\nu\mu}_{\rd(\nu,\lambda)-i-j}
\\
&=
\sum_i
\sum_{\mu,j}
(-1)^{\rd_I(\mu,\lambda)}\hat P^I_{\mu\lambda
j}t^{j}
Q^{\nu\mu}_{\rd(\nu,\lambda)-i-j}t^{\rd(\nu,\lambda)-i-j}
\\
&=
\sum_{\mu}
\sum_j
(-1)^{\rd_I(\mu,\lambda)}\hat P^I_{\mu\lambda
j}t^{j}
\sum_k
Q^{\nu\mu}_{k}t^{k}
\\
&=
\sum_{\mu}
(-1)^{\rd_I(\mu,\lambda)}\hat P^I_{\mu,\lambda
}(t^2)
Q^{\nu,\mu}(t^2),
\end{align*}
as desired.

\setcounter{equation}{0}
\noindent
(5.3)
Given a simple
$G_1T$-module, 
the formula (5.1.3)
is not necessarily 
accessible
to locate a simple factor
in the Loewy layers of $\hat\nabla_P(\hat L^P(\lambda))$. 
The following are particularly important factors in the study of the Frobenius direct image of
the structure sheaf of
$G/P$
\cite{HKR},
\cite{KY}.
Let $W^I=\{
w\in W\mid
\ell(ww')=\ell(w)+\ell(w')\ \forall w'\in W_I\}$, which forms a complete set of representatives
of
$W/W_I$.

\begin{prop}
Assume the Lusztig 
conjecture.
Let
$\lambda\in \Lambda$ be
$p$-regular.
If
$w\in W^I$,
$L((w\bullet\lambda)^0)\otimes
p(w^{-1}\bullet(w\bullet\lambda)^1)$
appears in the $(\ell(w)+1)$-st socle layer of
$\hat\nabla_P(\hat L^P(\lambda))$.

\end{prop}

\pf
In the commutative diagram (4.6.1) put
$\phi=\phi_{w^I}$. 
Write
$w^I=s_{i_1}s_{i_2}\dots
s_{i_m}$
in a reduced expression with
$m=\ell(w^I)$, and put
$y_r=s_{i_1}s_{i_2}\dots
s_{i_r}$ for
$r\leq m$.
Then
$y_r^{-1}w^I
=
s_{i_{r+1}}\dots
s_{i_m}\in W^I$.
Recall from \cite{AK89} that
$\phi_{w^I}:\hat\nabla_{w^I}((w^I\bullet
\lambda)\langle
w^I\rangle)\to
\hat\nabla(w^I\bullet
\lambda)$
is the composite
\begin{align*}
&
\hat\nabla_{w^I}((w^I\bullet
\lambda)\langle
w^I\rangle)
=
\hat\nabla_{s_{i_1}\dots
s_{i_m}}((w^I\bullet
\lambda)\langle
s_{i_1}\dots
s_{i_m}\rangle)
\overset{\phi_m}{\longrightarrow}
\\
&
\hat\nabla_{s_{i_1}\dots
s_{i_{m-1}}}((w^I\bullet
\lambda)\langle
s_{i_1}\dots
s_{i_{m-1}}\rangle)
\overset{\phi_{m-1}}{\longrightarrow}
\hat\nabla_{s_{i_1}\dots
s_{i_{m-2}}}((w^I\bullet
\lambda)\langle
s_{i_1}\dots
s_{i_{m-2}}\rangle)
\overset{\phi_{m-2}}{\longrightarrow}
\dots
\overset{\phi_{2}}{\longrightarrow}
\\
&
\hat\nabla_{s_{i_1}}((w^I\bullet
\lambda)\langle
s_{i_1}\rangle)
\overset{\phi_{1}}{\longrightarrow}
\hat\nabla(w^I\bullet
\lambda).
\end{align*}

Put $L=\soc
\hat\nabla_{y_r}((w^I\bullet
\lambda)\langle
y_r\rangle)\otimes-p(w^I\bullet0)$
and
$\phi_r'=\{(\phi_{r+1}\circ\dots\circ\phi_m)\otimes-p(w^I\bullet0)\}\vert_{^{w^I}\hat\nabla_P(\hat L^P(
\lambda))}$.
As $\phi_m'\ne0$ and as each
$\phi_i$ annihilates the socle
of its domain,
we must have
$\ell\ell(\im\phi_r')=\ell\ell(\hat\nabla_P(\hat L^P(
\lambda)))-(m-r)=
r+1$ by (4.6).
Then
$L
=\soc
(\im\phi'_r)=\rad_r(\im\phi_r')$,
which is a quotient of
$\rad_r{^{w^I}\hat\nabla}_P(\hat L^P(
\lambda))$.
Thus
$L$ lies in $\rad_r{^{w^I}\hat\nabla}_P(\hat L^P(
\lambda))$.
It follows from the rigidity of
$\hat\nabla_P(\hat L^P(
\lambda))$
that
$^{(w^I)^{-1}}\!L$ 
appears in its socle layer
of level
$\ell\ell(\hat\nabla_P(\hat L^P(
\lambda)))-r=\ell(w^I)+1-r=m+1-r=\ell(y_r^{-1}w^I)+1$.
Recall now from
\cite[1.2.4]{AK89} that
\begin{align*}
L
&=
\hat L((y_r^{-1}\bullet(w^I\bullet
\lambda))^0+p(y_r\bullet
(y_r^{-1}\bullet(w^I\bullet
\lambda))^1))\otimes-p(w^I\bullet0)
\\
&=
\hat L((y_r^{-1}w^I\bullet
\lambda)^0+p(y_r\bullet(y_r^{-1}w^I\bullet
\lambda)^1))\otimes-p(w^I\bullet0).
\end{align*}
Thus
\begin{align*}
^{(w^I)^{-1}}L
&=
L((y_r^{-1}w^I\bullet
\lambda)^0)\otimes
p\{
(w^I)^{-1}\{y_r\bullet(y_r^{-1}w^I\bullet
\lambda)^1-(w^I\bullet0)\}\}
\\
&=
L((y_r^{-1}w^I\bullet
\lambda)^0)\otimes
p\{(w^I)^{-1}y_r\bullet(y_r^{-1}w^I\bullet
\lambda)^1\}
\\
&=
L((y_r^{-1}w^I\bullet
\lambda)^0)\otimes
p\{(y_r^{-1}w^I)^{-1}\bullet(y_r^{-1}w^I\bullet
\lambda)^1\}.
\end{align*}

Finally,
we check that any $w\in W^I$ may be realized as $y_r^{-1}w^I$ as above.
Let $w\in W^I$.
As
$\ell(ww_I)=\ell(w)+\ell(w_I)$,
one can write
$w_0=s_{j_1}\dots
s_{j_r}ww_I$
with
$r=\ell(w_0)-\ell(w)-\ell(w_I)$.
Then
$w^I=w_0w_I=s_{j_1}\dots
s_{j_r}w$
with
$\ell(w^I)=r+\ell(w)$.
Thus, putting
$y_r=s_{j_1}\dots
s_{j_r}$ yields
$w=y_r^{-1}w^I$,
as desired.

\setcounter{equation}{0}
\noindent
(5.4)
{\bf Remark:}
This is a generalization of
\cite[1.5]{KY},
\cite{K09}
and
\cite[3.5]{K}.
In case $\lambda=0$
we constructed for 
$G$ of rank at most 2
\cite{KY}
or in case
$G=GL_{n+1}(\Bbbk)$ and $P$ is maximal parabolic such that
$G/P\simeq\bbP^n$ for any
$n\in\bbN$ \cite{K09}
a Karoubian complete strongly exceptional sequence 
$\{\cE_w\mid
w\in W^I\}$
for
the bounded derived category
of coherent
$\cO_{G_\bbC/P_\bbC}$-modules
out of
$G_1\Mod(L((w\bullet0)^0),
\soc_{\ell(w)+1}\hat\nabla_P(0))$,
where
$G_\bbC$
and
$P_\bbC$ are the groups over the complex number field
corresponding to
$G$ and $P$, respectively.
Our (5.3) assures at least that
$G_1\Mod(L((w\bullet0)^0),
\soc_{\ell(w)+1}\hat\nabla_P(0))\ne0$ in general
for large $p$.

\end{document}